\documentclass[preprint]{elsarticle}
\usepackage{lineno}

\journal{JOUARNAL's NAME}

\usepackage{slashbox}
\usepackage{amsmath}
\usepackage{amssymb}
\usepackage{amsthm}
\usepackage{graphicx}
\usepackage{epic,eepic,epsfig}
\usepackage{color}
\usepackage{subfigure} 
\usepackage{placeins}   
\usepackage{multirow}
\usepackage{epstopdf}
\usepackage{varwidth}
\usepackage{float}
\usepackage[table]{xcolor}
\usepackage{bm}
\usepackage[title]{appendix}
\usepackage{diagbox}
\usepackage{booktabs}
\usepackage[title]{appendix}

\newtheorem{theorem}{Theorem}[section]
\newtheorem{example}{Example}[section]
\newtheorem{lemma}{Lemma}[section]
\newtheorem{remark}{Remark}[section]

\newtheorem{assumption}{Assumption}[section]

\definecolor{tabclr}{cmyk}{0,0,1,0}

\allowdisplaybreaks[4]

\usepackage{caption2}
\usepackage{algorithm}
\usepackage{algorithmic}
\numberwithin{equation}{section}

\allowdisplaybreaks[4]

\makeatletter
\def\ps@pprintTitle{%
 \let\@oddhead\@empty
 \let\@evenhead\@empty
 \let\@oddfoot\@empty
 \let\@evenfoot\@empty
}
\makeatother

\begin{document}
\title{Exponential Runge--Kutta methods of collocation type for parabolic equations with time-dependent delay}

\author[bjut]{Qiumei Huang}
\ead{qmhuang@bjut.edu.cn}

\author[ins1,ins2]{Alexander Ostermann}
\ead{alexander.ostermann@uibk.ac.at}

\author[bjut,ins1]{Gangfan Zhong\corref{cor}}
\ead{gfzhong@emails.bjut.edu.cn}

\cortext[cor]{Corresponding author}
\address[bjut]{School of Mathematics, Statistics and Mechanics, Beijing University of Technology, Beijing 100124, China}

\address[ins1]{Department of Mathematics, University of Innsbruck, 6020 Innsbruck, Austria}
\address[ins2]{Digital Science Center, University of Innsbruck, 6020 Innsbruck, Austria}

\begin{abstract}
In this paper, exponential Runge--Kutta methods of collocation type  (ERKC) which were originally proposed in (Appl Numer Math 53:323--339, 2005) are extended to semilinear parabolic problems with time-dependent delay. Two classes of the ERKC methods are constructed and their convergence properties are analyzed. 
It is shown that methods with $s$ arbitrary nonconfluent collocation parameters achieve convergence of order $s$. 
Provided that the collocation para\-meters fulfill some additional conditions and the solutions of the problems exhibit sufficient temporal and spatial smoothness, we derive superconvergence results. Finally, some numerical experiments are presented to  illustrate our theoretical results.
\end{abstract}


\maketitle

\section{Introduction}\label{Sec:intro}
\setcounter{equation}{0}
In this study, we consider the following class of (abstract) semilinear parabolic problems with time-dependent delay
\begin{equation}\label{Eqn:problem}
\left\{\begin{aligned}
&u^\prime(t) + Au(t) = g(t,u(t),u^\tau(t)), && t > 0, \\
&u(t) = \phi(t),   && t\leq 0,
\end{aligned}\right.
\end{equation}
where $u^\tau(t)=u( t -\tau(t) )$. The delay $\tau:[0,+\infty)\to \mathbb{R}^+$ is assumed to satisfy the following conditions:
\begin{itemize}
\item[(i)] It is bounded away from 0, i.e., there exists a constant $\tau_0 > 0$ such that $\tau(t)\geq \tau_0$ for all $t\geq 0$.
\item[(ii)] The delayed argument $t - \tau(t)$ is strictly increasing for all $t\geq 0$.
\end{itemize}
This class of problems arises in many models such as in the animal movement model\-ling \cite{WangSalmaniw23:71}, in the dynamics of marine bacteriophage infections \cite{GourleyKuang04:550}, and in the control of diffusion processes \cite{Wang75:274}, among others. A systematic study of the existence and uniqueness of solutions to \eqref{Eqn:problem} can be found in \cite{Wu96Book}.

In the past, significant numerical research has been conducted for problems of the form \eqref{Eqn:problem} with constant delay; see, for example, \cite{HuangVandewalle12:579,vanPieter86:1,XuHuang23:57}. 
Compared to constant delay models, time-dependent delay models generally provide a more accurate description of the dynamic nature of real-world systems.
Never\-theless, the presence of time-dependent delays complicates both the mathematical modeling and the numeri\-cal solution, and relatively little numerical analysis has been performed for these models. Despite these difficulties, some studies have been conducted in this area. For example, Tang and Zhang \cite{TangZhang21:48,ZhangTang22:106233} constructed fully discrete $\theta$-methods to solve semilinear reaction-diffusion equations and Sobolev equations with time-dependent delay. However, their approach requires the absence of discontinuities in the true solution, which is a significant limitation, as discontinuities commonly arise in delay differential equations. Operator splitting for abstract delay equations has also been investigated in \cite{BatkaiCsomos13:315,BatkaiCsomos17:345,HansenStillfjord14:673}. 
Until now, there have been limited numerical studies on temporal high-order schemes for solving \eqref{Eqn:problem} that can account for discontinuities due to the presence of a time-dependent delay.

In recent years, exponential integrators have attracted considerable attention due to their effectiveness on stiff semilinear systems, such as \eqref{Eqn:problem} with $A$ having a large norm or being an unbounded operator, and the nonlinearity $g$ satisfying a local Lipschitz condition with a relatively small Lip\-schitz constant in a strip along the exact solution. 
The stiff convergence of exponential integrators, such as exponential Runge--Kutta methods \cite{HochbruckOstermann05:1069,HochbruckOstermann05:323,LuanOstermann13:3431}, exponential Rosenbrock methods \cite{HochbruckOstermann09:786} and exponential multistep methods \cite{HochbruckOstermann11:889}, has been investigated. 
Exponential integrators have also been extended to various areas, such as low-regularity integrators \cite{OstermannSchratz18:731},  exponential splitting methods \cite{HansenOstermann09:1485}, and phase field models \cite{HuangQiao24:116981}.
For a detailed overview of exponential integrators, we refer the reader to \cite{HochbruckOstermann10:209}. 
Some research on exponential integrators for (parabolic) differential equations with constant delay can be found in
\cite{DaiHuang23:350,Zhan21:113279,XuZhao10:2350,XuZhao11:1089,ZhanFang24:128978,ZhaoZhan16:96,ZhaoZhan20:366}. 
In this paper, we aim to extend the exponential Runge--Kutta methods of collocation type (ERKC) as presented in \cite{HochbruckOstermann05:323} to the larger class of problems \eqref{Eqn:problem}. For this purpose, we construct continuous extensions of the discrete solutions obtained by the ERKC methods, and we derive the convergence results. 
In contrast to the aforementioned works on exponential integrators for problems with constant delay, our research focuses on the following three specific topics.
\begin{itemize}
\item[(a)]
The existing works are restricted to nonlinearities of the form \( g:[0,T]\times X \times X\to X \), where \( X \) denotes the state space of problem~\eqref{Eqn:problem}, our analysis allows for a broader class of nonlinearities (see Assumption~\ref{Ass:V}).
\item[(b)]
In the constant-delay case, by choosing a uniform time step that divides \( \tau \), it follows that for every mesh point \( t_n \), the delayed argument \( t_n - \tau \) also coincides with a mesh point.
 This alignment simplifies the approximation of the delayed solution \( u(t_n - \tau) \) during time stepping. In contrast, for time-dependent delay, the delayed argument \( t_n - \tau(t_n) \) generally does not coincide with any mesh point. To address this, we propose two strategies for constructing continuous extensions of the numerical solution to approximate the delayed solution.
\item[(c)]
The existing studies are limited to uniform time stepping. However, for problems with time-dependent delay, variable step sizes are required to ensure that the discontinuity points of the solution are captured by the time mesh. The corresponding convergence analysis on nonuniform meshes is nontrivial and requires new techniques.
\end{itemize}

The outline of the paper is as follows. In Section~\ref{Sec:framework}, we summarize the employed abstract framework. Further, we construct two classes of $s$-stage ERKC methods for \eqref{Eqn:problem} in Section~\ref{Sec:ERKC}. The error analysis for these methods is carried out in Section~\ref{Sec:s}, and an error estimate of order $s$ is obtained. Provided that the collocation parameters satisfy some additional conditions and the problem possesses sufficient temporal and spatial regularity, we obtain superconvergence results of order $s+1$ and $s+1+\beta$ (see \eqref{Eqn:ass-s+2-1} for the definition of $\beta$) in Section~\ref{Sec:superconvergence}. Finally, we present some numerical experiments in Section~\ref{Sec:experiments} to illustrate the theoretical results.

\section{Analytical Framework}\label{Sec:framework}

Our analysis below will be based on an abstract formulation of \eqref{Eqn:problem} as an evolution equation with delay in a Banach space $(X,\|\cdot\|)$. Let $D(A)$ denote the domain of $A$ in $X$. Our basic assumptions on the operator $A$ are those of \cite{Henry81Book,Lunardi95Book}.

\begin{assumption}\label{Ass:sectorial}
Let $A: D(A) \rightarrow X$ be sectorial and $\overline{D(A)}=X$. Without restriction of generality, we assume that the spectrum $\sigma(A)$ of $A$, satisfies $\mathrm{Re}\,\sigma(A)>0$.
\end{assumption}

Under this assumption, the operator $-A$ is the infinitesimal generator of an analytic semigroup $\{\mathrm{e}^{-t A}\}_{t \geq 0}$, and the fractional powers of $A$ are well defined. We recall that $A$ satisfies the parabolic smoothing property  (see \cite{Henry81Book})
$$
\|\mathrm{e}^{-t A}\|_{X \leftarrow X}+\big\|t^\gamma A^\gamma \mathrm{e}^{-t A}\big\|_{X \leftarrow X}   \leq C, \quad  \gamma,t \geq 0,
$$
where for two Banach spaces $X$ and $Y$, the symbol $\|\cdot\|_{Y\gets X}$ denotes the operator norm from $X$ to $Y$.

For $0\leq \alpha <1$, we consider the Banach space $V=\{v\in X: A^\alpha v \in X\}
$ with norm $\|v\|_V=\| A^\alpha v\|$. Our basic assumption on $g$ is the following.
\begin{assumption}\label{Ass:V}
Let $g:[0, T] \times V \times V \rightarrow X$ be locally Lipschitz-continuous in a strip along the exact solution $u$. 
\end{assumption}

From Assumption \ref{Ass:V} we infer that there exists a real number $L(R, T)$ such that
$$
\|g(t, v_1,w_1)-g(t,v_2, w_2)\| \leq L \big( \|v_1-v_2\|_V+\|w_1-w_2\|_V\big)
$$
for all $t \in[0, T]$ and all $v_1,w_1,v_2,w_2\in V$ with 
$$
\max (\|v_1-u(t)\|_V,\| v_2  -u(t)\|_V,\| w_1 -u^\tau(t)\|_V,\|w_2-u^\tau(t)\|_V) \leq R.
$$

We provide an example fitting into this abstract framework.
\begin{example}
\rm Consider the following semilinear reaction-diffusion equation with time-dependent delay  in $(0,T]\times {\Omega}$, 
\begin{equation}\label{Eqn:problem-example}
\left\{\begin{aligned}
&\frac{\partial u}{\partial t}(t,{x})  -\Delta u(t,{x}) = \frac{1}{1+ u(t,{x})^2}+\frac{1}{1+ u(\frac{t}{2}-\frac{1}{2},{x})^2}, && t\in(0,T],~{x}\in\Omega, \\
&u(t,{x}) = 0,   && t\in(0,T],~{x}\in \partial \Omega, \\
&u(t,{x}) = \phi(t,{x}),   && t\in [-\tfrac{1}{2},0 ],~{x}\in\Omega,
\end{aligned}\right.
\end{equation}
where $\Delta$ is the Laplacian on the bounded open set $\Omega\subset \mathbb{R}^d$ with a sufficiently smooth boundary $\partial\Omega$. For an abstract formulation of \eqref{Eqn:problem-example}, we consider the linear (differential) operator on $D(A)\subset X$ defined by $Au=-\Delta u$ for all $u\in D(A)$. There are several choices of the space $X$, such as (see \cite{Lunardi95Book})
\begin{itemize}
\item[(i)] $X=L^2(\Omega)$, $D(A)= H_0^1(\Omega)\cap H^2(\Omega)$.
\item[(ii)] $X=C_0(\Omega)$, $D(A)= \{ u: u \in W^{2,p}(\Omega)~ \mbox{for all}~ p \geq 1,\,\Delta u\in C_0 (\Omega),\,u|_{\partial \Omega}=0 \}$.
\end{itemize}
In this way, \eqref{Eqn:problem-example} is reformulated as \eqref{Eqn:problem} with $g(t,v,w)=\frac{1}{1+v^2}+\frac{1}{1+ w^2}$, which satisfies the local Lipschitz condition in both cases (i) and (ii).
\end{example}

\section{ERKC methods}\label{Sec:ERKC}

An important feature of delay (partial) differential equations is that they possess solutions with discontinuities in the derivative at initial time and discontinuities in high order derivatives at subsequent points. 
This phenomenon has been well studied in the context of delay differential equations \cite{BellenZennaro13Book}. We provide a brief explanation of how this phenomenon also arises in delay parabolic equations. Suppose that $\phi\in C^2([-\tau_0,0];V)$ and $\tau\in C^2([0,T])$, and assume that the function $g$ is twice continuously Fr\'echet differentiable and that $\tfrac{\partial^2 g}{\partial {v^2}}(t,v,w)$ is locally Lipschitz continuous. Within the interval $[0,\xi_1]$, where $\xi_1-\tau(\xi_1)=0$, the initial value problem \eqref{Eqn:problem} reduces to a problem without delay
$$
u^\prime(t) + Au(t) = h(t,u(t)),\quad t \in [0,\xi_1], \quad u(0)=\phi(0),
$$
where $h(t,v)=g(t,v,\phi(t-\tau(t)))$. If $u(0)\in D(A)$ and the compatibility conditions 
$$
r_0=-Au(0)+h(0,u(0))\in D(A),\quad 
-Ar_0+\frac{\partial h}{\partial t}(0,u(0)) + \frac{\partial h}{\partial v}(0,u(0))r_0 \in V
$$
hold, we deduce from \cite[Theorem~7.1.2]{Lunardi95Book} that the solution satisfies $u\in C^2([0,\xi_1];V)$. Noting that the compatibility conditions also hold at $t=\xi_1$, it follows that $u\in C^2([\xi_1,\xi_2];V)$ as well, where $\xi_2-\tau(\xi_2)=\xi_1$. If $\phi^\prime(0^-)\neq u^\prime(0^+)$ in $V$, noting that $(t-\tau(t))^\prime\neq 0$ due to the strict monotonicity of $t-\tau(t)$, we obtain
$$
u^{\prime\prime}(\xi_1^-)-u^{\prime\prime}(\xi_1^+)
=\frac{\partial g}{\partial w}(\xi_1,u(\xi_1),\phi(\xi_1-\tau(\xi_1))) \big( \phi^\prime(0^-)-u^\prime(0^+)\big) \cdot (t-\tau(t))^\prime|_{t=\xi_1} \neq 0 ,
$$
where the notation $v(t^\pm)$ denotes the one-sided limits of $v$ as $t \to t^\pm$. The same argument can be extended to show that  $u:[\xi_1,T]\to V$ is twice continuously Fr\'echet differentiable. 
The above analysis indicates that discontinuities also arise in delay parabolic equations.

The set of potential discontinuity points $J= \{\xi_\mu:\mu \in \mathbb{N},~ \xi_\mu < T\}$ (which are called primary discontinuity points in the theory of delay differential equations \cite{BellenZennaro13Book}) are generated by the recursion 
$$
\xi_{\mu+1}-\tau(\xi_{\mu+1})=\xi_\mu,\quad \mu=0,1,\ldots;\quad \xi_0=0. 
$$
For $\xi_{m}<T \leq \xi_{m+1}$ with some $m \in \mathbb{N}$, we have $J=\{\xi_\mu\}_{\mu=0}^{m}$ and corresponding intervals $\mathcal{I}^\mu=[\xi_{\mu-1},\min\{\xi_\mu,T\}]$ for $\mu=1,2,\ldots,m+1$. At the points $\xi_\mu \in J$, the solutions of \eqref{Eqn:problem}, regardless of how regular the given functions $\phi$, $g$ and $\tau$ are, will generally exhibit a low degree of regularity.  In particular, if the functions \( \phi \), \( g \), and \( \tau \) are \( p \)-times continuously Fréchet differentiable and \( \frac{\partial^p g}{\partial v^p}(t,v,w) \) is locally Lipschitz continuous, then, provided that certain compatibility conditions are satisfied at $t = 0$, it follows from [20,~Theorem~7.1.2] that the solution \( u \in C^p(\mathcal{I}^\mu; V) \), $\mu=1,2,\ldots,m+1$.
Considering such discontinuity points, we choose the mesh
$$
I_h=\{t_n: 0 =t_0<t_1<\ldots<t_N=T\}
$$
such that $J\subseteq I_h$. Choosing nonconfluent collocation parameters $c_i \in [0,1]$ for $i=1,2\ldots,s$, we set 
$$
I_{n+1}=[t_{n},t_{n+1}],~~ h_{n+1}=t_{n+1}-t_n, ~~ t_{n,i}=t_{n}+c_ih_{n+1} ~~\mbox{and}~~ h =\max_{1\leq n \leq N} h_n.
$$

We now extend ERKC methods to \eqref{Eqn:problem}. The solution of \eqref{Eqn:problem} at time $t_{n}+h_{n+1}$ is given by the variation-of-constants formula
\begin{equation}\label{Eqn:variationofconstants}
 u(t_{n}+h_{n+1} =\mathrm{e}^{-h_{n+1}A}u(t_n)+\int_0^{h_{n+1}} \mathrm{e}^{-(h_{n+1}-\sigma)A} g(t_n+\sigma,u(t_n+\sigma),u^\tau(t_n+\sigma) )~\mathrm{d}\sigma.
\end{equation}
Following the idea in \cite{HochbruckOstermann05:323}, we approximate the function $g$ by the interpolation polynomial, which is 
$$
G_n(\sigma)= \sum_{j=1}^s \ell_j(\sigma/h_{n+1})G_{n,j},\quad\mbox{with}~
G_{n,j}=g(t_{n,j},U_{n,j},u^\tau(t_{n,j})),
$$
where $U_{n,j}$ is an approximation to $u(t_{n,j})$ which will be determined later, and $\ell_j$ is the Lagrange interpolation polynomial
$$
\ell_{j}(\rho) = \prod_{m=1,m\neq j}^{s} \frac{\rho - c_m}{c_j -c_m}.
$$
Substituting $g$ in \eqref{Eqn:variationofconstants} by $G_n$, results in the following approximation to $u(t_{n+1})$:
\begin{equation}\label{general-ERKC}
\begin{aligned}
U_{n+1}
& = \mathrm{e}^{-h_{n+1}A}U_n + h_{n+1} \sum_{i=1}^s b_i(-h_{n+1}A) G_{n,i},
\end{aligned}
\end{equation}
where $U_n\approx u(t_n)$ is given from the previous step and
\begin{equation}\label{general-ERKC-b}
\begin{aligned}
b_i(- h_{n+1} A) & = \int_0^{1} \mathrm{e}^{-(1-\xi)h_{n+1}A}   \ell_{i}(\xi) ~\mathrm{d}\xi \\
&=\frac{1}{ h_{n+1}} \int_0^{h_{n+1}} \mathrm{e}^{-(h_{n+1}-\sigma)A}   \ell_{i}(\sigma/h_{n+1}) ~\mathrm{d}\sigma.
\end{aligned}
\end{equation}
In a similar manner, $U_{n,i}\approx u(t_{n,i})$ is given by
\begin{equation}\label{general-ERKC-inner}
U_{n,i}
= \mathrm{e}^{-c_i h_{n+1}A}U_n + h_{n+1} \sum_{j=1}^s a_{ij}(- h_{n+1}A) G_{n,j},
\end{equation}
where
\begin{equation}\label{general-ERKC-a}
\begin{aligned}
a_{ij}(-h_{n+1} A) 
&  = \int_0^{c_i} \mathrm{e}^{-(c_i -\xi)h_{n+1}A}   \ell_{j}(\xi) ~\mathrm{d}\xi  \\
&=\frac{1}{ h_{n+1}} \int_0^{c_i h_{n+1}} \mathrm{e}^{-(c_i h_{n+1}-\sigma)A}   \ell_{j}(\sigma/h_{n+1}) ~\mathrm{d}\sigma.
\end{aligned}
\end{equation}

The scheme \eqref{general-ERKC}-\eqref{general-ERKC-a} is not computable yet, since the delayed values $u^\tau(t_{n,i})$ in $G_n$ still have to be approximated. One feasible approach is to construct a continuous approximation by using piecewise interpolation polynomials based on the known values $U_{k}$ and $U_{k,i}$ for $k=0,1,\ldots,n-1$ and $i=1,2,\ldots,s$. For this, we introduce an interpolation $\Pi_h U(t)$, which is a polynomial of degree $l$ when restricted to each interval $I_{k+1}$. The polynomial $\Pi_hU|_{I_{k+1}}$ is  based on the nodes of the set $S_{k+1}=\{t_k,t_{k+1}\}\cup\{t_{k,1},t_{k,2},\ldots,t_{k,s}\}$, satisfying
\begin{equation}\label{Eqn:PiU}
\Pi_h U (t_{k})=U_{k},~  \Pi_h U (t_{k+1})=U_{k+1}, ~  \mbox{and}
~ \Pi_h U (t_{k,i})=U_{k,i}~ ~\mbox{for}~1\leq i \leq s.
\end{equation}
Depending on the number of distinct interpolation points available, the degree $l$ is either $s-1$, $s$ or $s+1$. For $t\in [-\tau_0,0]$, we denote $\Pi_h U(t)=\phi(t)$. Therefore, for $t \in [t_n, t_{n+1}]$ and $\hat{t} = t - \tau(t)$, we have $\Pi_h U(\hat{t}) = \phi(\hat{t})$ if $\hat{t} \le 0$, or $\Pi_h U(\hat{t}) = (\Pi_h U)|_{I_{k+1}}(\hat{t})$ if $\hat{t} \in I_{k+1}$ for some $0 \le k \le n-1$.
Replacing the evaluations $u^\tau(t_{n,i})$ in $G_{n}$ by $(\Pi_h U)^\tau(t_{n,i})$, we arrive at the scheme
\begin{equation}\label{ERKC-1x}
\begin{aligned}
U_{n+1}
& = \mathrm{e}^{-h_{n+1}A}U_n + h_{n+1} \sum_{i=1}^s b_i(-h_{n+1}A) G^{\rm I}_{n,i}, \\
U_{n,i} &
= \mathrm{e}^{-c_i h_{n+1}A}U_n + h_{n+1} \sum_{j=1}^s a_{ij}(- h_{n+1}A) G^{\rm I}_{n,j},\quad 1\leq i \leq s, \\
G^{\rm I}_{n,j}  & =g(t_{n,j},U_{n,j},(\Pi_hU)^\tau(t_{n,j})),\quad 1\leq j \leq s,
\end{aligned}
\end{equation}
which we will call ERKC-I method, as it employs an interpolation polynomial to approximate the delayed solution. For $t \geq 0$, $\Pi_h U(t)$ is a continuous extension of the discrete numerical solution obtained in \eqref{ERKC-1x}.

The scheme \eqref{ERKC-1x} is one option among several possible schemes for extending ERKC methods to \eqref{Eqn:problem}.
We shall introduce an alternative approach for constructing an approximation to the delayed solution $u^\tau(t_{n,i})$ in $G_n$. To distinguish it from the previous introduced numerical solution $U_{n+1}$ and $U_{n,i}$ in \eqref{ERKC-1x}, we denote the approximate solutions of the new method by $W_{n+1}$ and $W_{n,i}$:
$$
\begin{aligned}
W_{n+1}
& = \mathrm{e}^{-h_{n+1}A}W_n + h_{n+1} \sum_{i=1}^s b_i(-h_{n+1}A)g(t_{n,i},W_{n,i},W^\tau(t_{n,i})), \\
W_{n,i} &
= \mathrm{e}^{-c_i h_{n+1}A}W_n + h_{n+1} \sum_{j=1}^s a_{ij}(- h_{n+1}A) g(t_{n,j},W_{n,j},W^\tau(t_{n,j})),~~ 1\leq i \leq s,
\end{aligned}
$$
with $W(t)$ still to be defined. Instead of using piecewise interpolation polynomials, we construct a continuous solution $W(t)$ by setting $W(t)= \phi(t)$ for $t\in [-\tau_0,0]$ and, for $t=t_{k} + \theta h_{k+1}$ with some $k\in\{1,2,\ldots,n\}$ and $\theta\in(0,1]$,
\begin{equation}\label{Eqn:VV}
W(t_k+\theta h_{k+1})=\mathrm{e}^{-\theta h_{k+1}A}W_k + h_{k+1} \sum_{i=1}^s b_i(\theta;- h_{k+1}A) g(t_{k,i},W_{k,i},W^\tau(t_{k,i})),
\end{equation}
where
\begin{equation}\label{general-ERKC-btheta}
\begin{aligned}
b_i(\theta;- h_{k+1} A) 
&  =   \int_{0}^{\theta} \mathrm{e}^{-(\theta-\xi)h_{k+1}A}\ell_{i}( \xi)~\mathrm{d}\xi   \\
&=\frac{1}{ h_{k+1}} \int_0^{\theta h_{k+1}} \mathrm{e}^{-(\theta h_{k+1}-\sigma)A}   \ell_{i}(\sigma/h_{k+1}) ~\mathrm{d}\sigma.
\end{aligned}
\end{equation}
We thus arrive at the scheme
\begin{equation}\label{ERKC-2x}
\begin{aligned}
W_{n+1}
& = \mathrm{e}^{-h_{n+1}A}W_n + h_{n+1} \sum_{i=1}^s b_i(-h_{n+1}A) G^{\rm C}_{n,i}, \\
W_{n,i} &
= \mathrm{e}^{-c_i h_{n+1}A}W_n + h_{n+1} \sum_{j=1}^s a_{ij}(- h_{n+1}A) G^{\rm C}_{n,j},\quad 1\leq i \leq s, \\
G^{\rm C}_{n,j} & =g(t_{n,j},W_{n,j},W^\tau(t_{n,j})),\quad 1\leq j \leq s.
\end{aligned}
\end{equation}
For $A=0$, this scheme reduces to a collocation method, in which the delayed solution is approx\-i\-mated by the collocation polynomials. For this reason, we refer to it as the ERKC-C method.
Noting that 
\begin{equation}\label{Eqn:ab}
b_i(-h_{n+1}A)=b_i(1;-h_{n+1}A)\quad\mbox{and}\quad  a_{ij}(-h_{n+1}A)=b_{j}(c_i;-h_{n+1}A),
\end{equation}
we immediately obtain $W(t_{k+1})=W_{k+1}$ and $W(t_{k,i})=W_{k,i}$ for $k=0,1,\ldots,N-1$ and $i=1,2,\ldots,s$. For each fixed $\theta \in [0,1]$, the coefficients $b_i(\theta;- h_{n+1}A)$, as defined in \eqref{general-ERKC-btheta}, are linear combinations of the functions $\varphi_j(-\theta h_{n+1} A)$, where
\begin{equation}\label{Eqn:varphi}
\varphi_j(-t A)  =\int_0^1\mathrm{e}^{-(1-\xi)tA}\frac{\xi^{j-1}}{(j-1)!}
~\mathrm{d}\xi=\frac{1}{t^j} \int_0^t \mathrm{e}^{-(t-\sigma) A} \frac{\sigma^{j-1}}{(j-1)!} ~\mathrm{d} \sigma, \quad 1 \leq j \leq s.
\end{equation}
The coefficients in the linear combination depend on $\theta$; see, e.g., Example \ref{Exa:s=2}. 
Recalling the parabolic smoothing property of $A$, we have the estimate
\begin{equation}\label{Eqn:smoothing-btheta}
\| t^\gamma A^\gamma b_i(\theta;- t A) \|_{X\leftarrow X} \leq C,\quad 0\leq  \gamma \leq 1,~t\geq0.
\end{equation}

\begin{example}\rm
For $s=1$, the weight \eqref{general-ERKC-btheta} is given by
$$
b_1(\theta;z) = \theta \varphi_1(\theta z).
$$
The choice $c_1=1$ leads to the exponential backward Euler method. The ERKC-I form of this scheme is 
$$
U_{n+1}
 = \mathrm{e}^{-h_{n+1}A}U_n + h_{n+1}\varphi_1(-h_{n+1}A) g(t_{n+1},U_{n+1},(\Pi_hU)^\tau(t_{n+1})),  
 $$
where $\Pi_hU(t)=(1-\theta)U_k + \theta U_{k+1}$, for $t=t_k+\theta h_{k+1} \in I_{k+1}$, $0\leq k \leq n-1$.
The ERKC-C variant is
$$
W_{n+1}
 = \mathrm{e}^{-h_{n+1}A}W_n + h_{n+1}\varphi_1(-h_{n+1}A) g(t_{n+1},W_{n+1},W^\tau(t_{n+1})), \\
$$
where $W(t)=\mathrm{e}^{-\theta h_{k+1}A}W_k + \theta h_{k+1} \varphi_1(-\theta h_{n+1}A) G^{\rm C}_{k,1}$, for $t=t_k+\theta h_{k+1} \in I_{k+1}$, $0\leq k \leq n-1$.
\end{example}

\begin{example}\label{Exa:s=2}\rm
For $s=2$, the weights \eqref{general-ERKC-btheta} are given by
$$
\begin{aligned}
b_1(\theta;z) & =\frac{c_2}{c_2-c_1} \cdot \theta \varphi_1(\theta z)-\frac{1}{c_2-c_1} \cdot \theta^2 \varphi_2(\theta z), \\
b_2(\theta;z) & =-\frac{c_1}{c_2-c_1} \cdot \theta \varphi_1(\theta z)+\frac{1}{c_2-c_1}\cdot \theta^2 \varphi_2(\theta z).
\end{aligned}
$$
The choice $c_1=1/3$, $c_2=1$ leads to the exponential Radau IIA method. The ERKC-I form of this scheme is 
$$
\begin{aligned}
U_{n+1}
& = \mathrm{e}^{-h_{n+1}A}U_n \!\!\!\!\!\!\!\!\!\! && \!\!\!\!\!  \!\!\! +h_{n+1}b_1(1;-h_{n+1}A)g(t_{n,1},U_{n,1},(\Pi_hU)^\tau(t_{n,1})) \\
& \phantom{=\mathrm{e}^{-h_{n+1}A}U_n} && \!\!\!\!\!  \!\!\! + h_{n+1}b_2(1;-h_{n+1}A)g(t_{n+1},U_{n+1},(\Pi_hU)^\tau(t_{n+1})), \\
U_{n,1}
& = \mathrm{e}^{-\tfrac{1}{3}h_{n+1}A}U_n \!\!\!\!\!\!\!\!\!\! && \!\!\!\!
+h_{n+1}b_1(\tfrac{1}{3};-h_{n+1}A)g(t_{n,1},U_{n,1},(\Pi_hU)^\tau(t_{n,1})) \\
& \phantom{= \mathrm{e}^{-\tfrac{1}{3}h_{n+1}A}U_n} && \!\!\!\! + h_{n+1}b_2(\tfrac{1}{3};-h_{n+1}A)g(t_{n+1},U_{n+1},(\Pi_hU)^\tau(t_{n+1})), 
\end{aligned}
$$
where $\Pi_hU(t)= 3 (\theta-\tfrac{1}{3} )(\theta-1)U_{k}
 -\tfrac{9}{2}\theta(\theta-1)U_{k,1} + \tfrac{3}{2}\theta (\theta-\tfrac{1}{3} )U_{k+1}$, for $t=t_k+\theta h_{k+1} \in I_{k+1},~ 0\leq k \leq n-1$. The ERKC-C variant is
$$
\begin{aligned}
W_{n+1}
& = \mathrm{e}^{-h_{n+1}A}W_n \!\!\!\!\!\!\!\!\!\! && \!\!\!\!\!  \!\!\! +h_{n+1}b_1(1;-h_{n+1}A)g(t_{n,1},W_{n,1},W^\tau(t_{n,1})) \\
& \phantom{=\mathrm{e}^{-h_{n+1}A}V_n} && \!\!\!\!\!  \!\!\! + h_{n+1}b_2(1;-h_{n+1}A)g(t_{n+1},W_{n+1},W^\tau(t_{n+1})), \\
W_{n,1}
& = \mathrm{e}^{-\tfrac{1}{3}h_{n+1}A}W_n \!\!\!\!\!\!\!\!\!\! && \!\!\!\!\! +h_{n+1}b_1(\tfrac{1}{3};-h_{n+1}A)g(t_{n,1},W_{n,1},W^\tau(t_{n,1})) \\
& \phantom{= \mathrm{e}^{-\tfrac{1}{3}h_{n+1}A}U_n} && \!\!\!\!\! + h_{n+1}b_2(\tfrac{1}{3};-h_{n+1}A)g(t_{n+1},W_{n+1},W^\tau(t_{n+1})), 
\end{aligned}
$$
where 
$
W(t)=\mathrm{e}^{-\theta h_{k+1}A}W_k + h_{k+1} b_1(\theta;-h_{k+1}A) G_{k,1}^{\rm C}  +h_{k+1} b_2(\theta;-h_{k+1}A) G_{k,2}^{\rm C}, 
$
for $t=t_k+\theta h_{k+1} \in I_{k+1}$, $0\leq k \leq n-1$.
\end{example}

 \section{Error analysis}\label{Sec:s}
Based on the construction of the stages $U_{n,i}$ of the ERKC-I method, we construct the following auxiliary function on $I_{n+1}$: 
\begin{equation}\label{ERKC-1x-aux}
U(t_n+\theta h_{n+1}) 
= \mathrm{e}^{-\theta h_{n+1}A}U_n + h_{n+1} \sum_{i=1}^s b_i(\theta; -h_{n+1}A) G^{\mathrm{I}}_{n,i}, ~~ 0\leq \theta \leq 1.
\end{equation}
Due to \eqref{Eqn:ab}, $U(t)$ satisfies $U(t_{n,i})=U_{n,i}$ for $1\leq i \leq s$ and $U(t_{n+1})=U_{n+1}$, which indicates that its interpolation $\Pi_h U(t)$ is the same polynomial as defined in Section~\ref{Sec:ERKC}.
This function is very convenient for the numerical analysis.

Henceforth, we suppose that \eqref{Eqn:problem} admits a solution $u:[0,T]\to V$ that is sufficiently smooth in each segment $\mathcal{I}^\mu$, with its derivatives belonging to $V$.
Assuming that the nonlinearity $g(t, v,w)$ is sufficiently often Fr\'echet differentiable in a strip along the exact solution, that the delay $\tau(t)$ is sufficiently smooth in $[0,T]$ and that $\phi(t)$ is sufficiently smooth in $[-\tau_0,0]$, it follows that the composition
\begin{equation}\label{Eqn:definition-f}
f:[0,T]\to X : t \mapsto f(t)=g(t,u(t),u^\tau(t)) 
\end{equation}
is sufficiently smooth in each segment $\mathcal{I}^\mu$ as well.

Inserting the exact solution into \eqref{ERKC-1x-aux} shows
\begin{align}
u(t_n+\theta h_{n+1})& = \mathrm{e}^{-\theta h_{n+1}A}u(t_n) +  h_{n+1}\sum_{i=1}^s b_i(\theta;- h_{n+1}A)f(t_n+c_ih_{n+1}) + \delta_{n+1}(\theta), \label{Eqn:exact-to-1x}
\end{align}
with defects $\delta_{n+1}(\theta)$. To derive the explicit expression for the defects, we substitute the Taylor series of $f$ given by
$$
f(t_n+\sigma)=\sum_{j=1}^s \frac{\sigma^{j-1}}{(j-1)!} f^{(j-1)}(t_n^+)+\int_0^\sigma \frac{(\sigma-\rho)^{s-1}}{(s-1)!} f^{(s)}(t_n+\rho) ~ \mathrm{d} \rho 
$$
into \eqref{Eqn:exact-to-1x} and into the exact solution (represented by the variation-of-constants formula)
\begin{equation*}
u(t_{n}+\theta h_{n+1})=\mathrm{e}^{-\theta h_{n+1}A}u(t_n)+\int_0^{\theta h_{n+1}} \mathrm{e}^{-(\theta h_{n+1}-\sigma)A} f(t_n+\sigma) ~\mathrm{d}\sigma.
\end{equation*}
Thus, we obtain
\begin{align}
 u(t_n+\theta h_{n+1}) 
&  = \mathrm{e}^{-\theta h_{n+1}A}u(t_n) +  h_{n+1}\sum_{i=1}^s b_i(\theta;- h_{n+1}A)  \sum_{j=1}^s \frac{(c_i h_{n+1})^{j-1}}{(j-1)!} f^{(j-1)}(t_n^+) \notag  \\
& \quad +  h_{n+1}\sum_{i=1}^s b_i(\theta;- h_{n+1}A)  \int_0^{c_i h_{n+1}} \frac{(c_i h_{n+1}-\rho)^{s-1}}{(s-1)!} f^{(s)}(t_n+\rho) ~ \mathrm{d} \rho  + \delta_{n+1}(\theta), \label{Eqn:pre-defect-1}
\end{align}
and
\begin{align}
u(t_n+\theta h_{n+1})&=  \mathrm{e}^{-\theta h_{n+1} A} u(t_n)+\sum_{j=1}^{s}(\theta h_{n+1})^j \varphi_j(-\theta h_{n+1} A) f^{(j-1)}(t_n^+) \notag \\
&\quad +\int_0^{\theta h_{n+1}} \mathrm{e}^{-(\theta h_{n+1}-\sigma) A} \int_0^\sigma \frac{(\sigma-\rho)^{s-1}}{(s-1)!} f^{(s)}(t_n+\rho) ~\mathrm{d} \rho ~\mathrm{d} \sigma, \label{Eqn:pre-defect-2}
\end{align}
respectively. Subtracting \eqref{Eqn:pre-defect-2} from \eqref{Eqn:pre-defect-1} and using
\begin{equation}\label{Eqn:order-conditions}
 \sum_{i=1}^s b_i(\theta; -h_{n+1}A)  \frac{c_i^{j-1}}{(j-1)!}= \theta^{j}\varphi_{j}(-\theta h_{n+1}A),\quad 1\leq j \leq s
\end{equation}
(which is an immediate consequence of \eqref{general-ERKC-btheta} and \eqref{Eqn:varphi}), we have
\begin{align}
\delta_{n+1}(\theta)= & \int_0^{\theta h_{n+1}} \mathrm{e}^{-(\theta h_{n+1}-\sigma) A} \int_0^\sigma \frac{(\sigma-\rho )^{s-1}}{(s-1)!} f^{(s)}(t_n+ \rho ) ~\mathrm{d} \rho ~\mathrm{d} \sigma  \notag \\
& - h_{n+1} \sum_{i=1}^s b_i(\theta;- h_{n+1} A) \int_0^{c_i h_{n+1}} \frac{(c_i h_{n+1}-\sigma)^{s-1}}{(s-1)!} f^{(s)}(t_n+\sigma) ~\mathrm{d} \sigma. \label{Eqn:defect}
\end{align}
We note that the identities \eqref{Eqn:order-conditions} are usually called order conditions. Using \eqref{Eqn:smoothing-btheta}, the following estmate of $\delta_{n+1}(\theta)$ is obtained
\begin{equation}\label{Eqn:defect-estimate}  
h_{n+1}^{\alpha} \max_{0\leq \theta \leq 1} \| \delta_{n+1}(\theta) \|_V + \max_{0\leq \theta \leq 1} \| \delta_{n+1}(\theta)\|_X
\leq C h_{n+1}^s\int_{t_n}^{t_{n+1}} \|f^{(s)}(\sigma)\|_X ~\mathrm{d}\sigma.
\end{equation}
Let $e(t)=U(t)-u(t)$ denote the difference between the auxiliary function \eqref{ERKC-1x-aux} and the exact solution. Subtracting \eqref{Eqn:exact-to-1x} from \eqref{ERKC-1x-aux} gives the error recursion
\begin{equation}\label{Eqn:main-error-1x}
  e(t_n+\theta h_{n+1}) 
  = \mathrm{e}^{-\theta h_{n+1}A}e(t_n) +  h_{n+1} \sum_{i=1}^s b_i(\theta;- h_{n+1}A)\Big( g(t_{n,i},U_{n,i}, (\Pi_h{U})^{\tau}_{n,i} )  - f(t_{n,i}) \Big) -\delta_{n+1}(\theta),  
\end{equation}
where $(\Pi_hU)^{\tau}_{n,i}=\Pi_hU(t_{n,i} - \tau(t_{n,i}))$. 
Similar to the construction of $\Pi_h U(t)$, we introduce $\Pi_h u(t)$, which is an interpolation polynomial of $u(t)$ when restricted to the interval $I_{k+1}$ ($k=0,1,\ldots,N-1$). This polynomial is also based on the nodes of the set $S_{k+1}=\{t_k,t_{k+1}\}\cup\{t_{k,1},t_{k,2},\ldots,t_{k,s}\}$, satisfying
\begin{equation}\label{Eqn:Piu}
\Pi_h u (t_{k})=u(t_{k}),~ \Pi_h u (t_{k+1})=u(t_{k+1}), ~ \mbox{and}
~ \Pi_h u (t_{k,i})=u(t_{k,i})~~\mbox{for}~1\leq i \leq s.
\end{equation}
The error of the numerical solution $\Pi_hU(t)$ in each interval $I_{k+1}$ can be controlled by
\begin{align}
 \max_{t\in I_{k+1}} \| \Pi_hU (t)-u(t)\|_V 
&   \leq \max_{t\in I_{k+1}} \|\Pi_hU(t) - \Pi_hu(t) \|_V  +\max_{t\in I_{k+1}}\| \Pi_hu(t) - u(t)\|_V \notag \\
&   \leq C \max_{\sigma \in S_{k+1}} \|U(\sigma)-u(\sigma)\|_V+Ch^{s}_{k+1} \sup_{\sigma\in (t_{k},t_{k+1})}\|u^{(s)}(\sigma)\|_V , \label{Eqn:interpolation-s}
\end{align}
provided that $u^{(s)}\in L^\infty(0,T;V)$, where the constant $C$ is independent of $I_{k+1}$ and the step size sequence.

The convergence result for the ERKC-I methods is presented in the next theorem, under the following assumption on the step size sequence: there is a constant $\rho>1$ such that
\begin{equation}\label{Eqn:step-size}
h_j\leq \varrho h_{j+1} \quad~\mbox{for all}~j.
\end{equation}
This is only a weak constraint on the step size sequence. Note that new step size $h_{j+1}$ is only constrained from above by the maximum step size $h$ and from below by the ratio $h_{j}/\rho$. For example, both quasi-uniform meshes and graded meshes satisfy this requirement.

\begin{theorem}\label{Thm:1x-s-global-error}
Let the initial value problem \eqref{Eqn:problem} satisfy Assumptions \ref{Ass:sectorial}-\ref{Ass:V}, and consider for its numerical solution the ERKC-I method \eqref{ERKC-1x}. Let the step size sequence $\{h_j\}_{j=1}^N$ satisfy the condition \eqref{Eqn:step-size}.  
If $u:\mathcal{I}^\mu\to V$ and $f:\mathcal{I}^\mu\to X$ \eqref{Eqn:definition-f} are both $s$ times differentiable for $\mu=1,2,\ldots,m+1$ with $u^{(s)}\in L^\infty(0,T;V)$, $f^{(s)}\in L^\infty(0,T;X)$, then for $h = \max_{1\leq j \leq N} h_j$ sufficiently small, the numerical solution \eqref{Eqn:PiU} is determined uniquely and the error of $\Pi_hU(t)$ satisfies 
\begin{equation*}
\max_{0\leq t \leq T}\|\Pi_h U(t)-u(t)\|_V \leq    Ch^s \bigg( \sup_{\sigma\in[0,T]}\|f^{(s)}(\sigma)\|_X +  \sup_{\sigma\in[0,T]} \|u^{(s)}(\sigma)\|_V\bigg) ,
\end{equation*}
where the constant $C$  depends on $T$, but is independent of the step size sequence.
\end{theorem}
\begin{proof}
For $n=0,1,\ldots,N-1$, 
the error in the interval $I_{n+1}$ fulfills the preliminary estimate
\begin{align}
 \|e(t_n+\theta h_{n+1})\|_V 
& \leq C\|e(t_n)\|_V +  CL\cdot h_{n+1}^{1-\alpha} \sum_{j=1}^s\bigg( \|e(t_{n,j})\|_V   + \|(\Pi_hU)^\tau(t_{n,j})-u^\tau(t_{n,j})\|_V\bigg) \notag \\
&\quad +
\|\delta_{n+1}(\theta)\|_V,\quad \theta\in[0,1].    \label{Eqn:error-in-V} 
\end{align}
Thus, for $h$ sufficiently small, we obtain
\begin{align}
 \|e(t_{n,i})\|_V
& \leq C\bigg(\|e(t_n)\|_V + \sum_{j=1}^s   \|\delta_{n+1}(c_j)\|_V\bigg)
+\sum_{j=1}^s \|(\Pi_hU)^\tau(t_{n,j})-u^\tau(t_{n,j})\|_V. \label{Eqn:inner-error-in-V} 
\end{align}
The proof proceeds by induction on the intervals $[0,\xi_{\mu}]$, $\mu=1,2,\ldots,m+1$. We first consider the interval $[0,\xi_1]$ in which case the problem reduces to a parabolic problem without delay.
When the mesh point $t_n$ belongs to $[0,\xi_1)$, we arrive at
\begin{equation}\label{Eqn:s-01-inner}
 \|e(t_{n,i})\|_V
  \leq C\bigg(\|e(t_n)\|_V + \sum_{j=1}^s   \|\delta_{n+1}(c_j)\|_V\bigg).
\end{equation}
Solving the error recursion \eqref{Eqn:main-error-1x} gives
$$
\begin{aligned}
 e(t_n) 
&=  \mathrm{e}^{-(t_n-t_0)A}e(t_0)   +  \sum_{j=1}^{n}\mathrm{e}^{-(t_n-t_j)A}\bigg(  \sum_{i=1}^s h_jb_i(-h_jA) \Big(g(t_{j-1,i},U_{j-1,i},(\Pi_hU)^\tau_{j-1,i}) 
- f(t_{j-1,i}) \Big) \\
&\quad -\delta_{j}(1)\bigg).
\end{aligned}
$$
Recalling the parabolic smoothing property of $A$ and using $e(t_0)=0$, we estimate this term in $V$ by
\begin{equation}\label{Eqn:recursion-in-V}
\begin{aligned}
\|e(t_n)\|_V & \leq  CL\sum_{j=1}^{n-1} (t_n-t_j)^{-\alpha} h_j \sum_{i=1}^s \big( 
\|e(t_{j-1,i})\|_V + \| (\Pi_hU)^\tau_{j-1,i}-u^\tau_{j-1,i}\|_V \big)     \\
& \quad  + CL\cdot h_{n}^{1-\alpha} \sum_{i=1}^s\big( \|e(t_{n-1,i})\|_V +\| (\Pi_hU)^\tau_{n-1,i}-u^\tau_{n-1,i}\|_V \big)    \\
  &\quad + C\sum_{j=1}^{n-1} (t_n-t_j)^{-\alpha} \|\delta_j(1)\|_X
  +    \|\delta_n(1)\|_V  . 
\end{aligned}
\end{equation}
Note that $h_j\leq \varrho h_{j+1}$ and $(t_n-t_j)^{-\alpha} \leq (1+\varrho)^\alpha(t_n-t_{j-1})^{-\alpha}$.
Using \eqref{Eqn:defect-estimate}, \eqref{Eqn:s-01-inner} and the fact that $\Pi_h U(t)=u(t)$ for $t\in [-\tau_0,0]$, we get by applying the discrete Gronwall lemma (Lemma \ref{Lem:gronwall}) the following result: for any mesh point $t_n$ belonging to $[0,\xi_{1}]$, it holds
\begin{equation}\label{Eqn:error-at-node-s-01} 
\|e(t_n)\|_{V}  \leq Ch^{s} \sup_{\sigma\in[0,\xi_{1}]} \|f^{(s)}(\sigma)\|_X.
\end{equation}
Combining \eqref{Eqn:error-in-V}, \eqref{Eqn:s-01-inner} and \eqref{Eqn:error-at-node-s-01}, we obtain
\begin{equation}\label{Eqn:error-global-s-01} 
\max_{t\in[0,\xi_1]} \|U(t)-u(t)\|_{V}  \leq Ch^{s} \sup_{\sigma\in[0,\xi_{1}]} \|f^{(s)}(\sigma)\|_X.
\end{equation}

Now, let
$$
B_{\mu_1,\mu_2}= \sup_{\sigma\in[0,\xi_{\mu_1}]} \|f^{(s)}(\sigma)\|_X  +   \sup_{\sigma\in[0,\xi_{\mu_2}]} \|u^{(s)}(\sigma)\|_V .
$$
In order to carry out the induction step, we show that if
\begin{equation}
\begin{aligned}
\max_{t\in[0,\xi_l]} \|U(t)-u(t)\|_V & \leq Ch^{s}B_{l,l-1}
\label{Eqn:s-assume-0k} 
\end{aligned}
\end{equation}
holds for $l=\mu$, then it also holds for $l=\mu+1$. Note that a combination of \eqref{Eqn:interpolation-s} and \eqref{Eqn:s-assume-0k} leads to
\begin{equation}
\begin{aligned}
\max_{t\in[0,   \xi_\mu]}\| \Pi_hU (t)-u(t)\|_V &  \leq Ch^{s}B_{\mu,\mu}.
\label{Eqn:interpolation-s+1-xi0xik}
\end{aligned}
\end{equation}
For a mesh point $t_{n}$ belonging to $[\xi_\mu,\xi_{\mu+1})$, from \eqref{Eqn:inner-error-in-V} and \eqref{Eqn:interpolation-s+1-xi0xik}, we obtain
\begin{equation}\label{Eqn:s-0k+1-inner}
\begin{aligned}
 \|e(t_{n,i})\|_V &\leq C\|e(t_n)\|_V +Ch^{s}{B_{\mu+1,\mu}}.
\end{aligned}
\end{equation}
Solving the error recursion \eqref{Eqn:main-error-1x} for $t_n\in  (\xi_\mu,\xi_{\mu+1})$ gives
$$
\begin{aligned}
 e(t_n) 
&= \mathrm{e}^{-(t_n-\xi_\mu)A}e(\xi_\mu) + \sum_{j=\chi(\xi_{\mu})+1}^{n}\mathrm{e}^{-(t_n-t_j)A}\bigg(  \sum_{i=1}^s h_jb_i(-h_jA) \Big(g(t_{j-1,i},U_{j-1,i},(\Pi_hU)^\tau_{j-1,i}) \\
&\qquad  - f(t_{j-1,i}) \Big) -\delta_{j}(1)\bigg),
\end{aligned}
$$
where $t_{\chi(\xi_\mu)}=\xi_\mu$.
Now, estimating this term in $V$ similarly to \eqref{Eqn:recursion-in-V}, using \eqref{Eqn:interpolation-s+1-xi0xik} and \eqref{Eqn:s-0k+1-inner} and applying the discrete Gronwall lemma (Lemma \ref{Lem:gronwall}) leads to the desired result: for a mesh point $t_n$ belonging to $[\xi_\mu,\xi_{\mu+1}]$,
\begin{equation}\label{Eqn:error-at-node-s-0k+1} 
\|e(t_n)\|_{V}  \leq  Ch^{s}B_{\mu+1,\mu}.
\end{equation}
Combining \eqref{Eqn:error-in-V} and \eqref{Eqn:interpolation-s+1-xi0xik}-\eqref{Eqn:error-at-node-s-0k+1}, we finally obtain
\begin{equation}\label{Eqn:error-global-s-0k+1} 
\max_{t\in[0,\xi_{  \mu+1}]} \|U(t)-u(t)\|_{V}  \leq Ch^{s} B_{\mu+1,\mu},
\end{equation}
which completes the induction.

The existence and uniqueness of the numerical solutions $U_{n,i}~(n=0,1,\ldots,N-1,~i=1,2,\ldots,s)$ can be verified by fixed-point iteration (cf.~\cite{OstermannThalhammer02:367}), proceeding by induction on the intervals $\mathcal{I}^\mu$. Since $g$ is locally Lipschitz continuous, the convergence proof proceeds as long as the numerical solution remains within a prescribed neighborhood of the exact solution, which can be ensured by choosing the time step size $h$ sufficiently small.
\end{proof}

The above proof relies on the following discrete Gronwall-type inequality involving weakly singular kernels and non-uniform time discretization, which is a simple extension of Theorem 6.1 in~\cite{DixonMcKee86:535}.

\begin{lemma}\label{Lem:gronwall}
Assume the sequences $\{h_j\}_{j=1}^N$, $\{a_{n,j}\}_{j=1}^{n-1}$ $(n\leq N)$ and $\{b_{j}\}_{j=1}^{N}$ are nonnegative and satisfy
$$
\begin{aligned}
 a_{n,j}  \leq a(t_n-t_j)^\gamma  \quad \mbox{and}\quad
b_j   \leq b  ,
\end{aligned}
$$
where $a>0$, $b>0$ and $-1<\gamma \leq 0$. Further assume that the nonegative sequence $\{\varepsilon_j\}_{j=1}^N$ satisfies the inequalities
$$
\varepsilon_n \leq   \sum_{j=1}^{n-1} h_{j+1} a_{n,j} \varepsilon_j + b_n,\quad\mbox{for}~ 1\leq n \leq N.
$$
Then the following estimate holds
$$
\varepsilon_n \leq Cb ,\quad \mbox{for}~1\leq n \leq N,
$$
where the constant $C$ depends on $a$, $\gamma$ and on $t_N$.
\end{lemma}

The convergence result for ERKC-C methods is presented in the next theo\-rem. The construction of $U(t)$ \eqref{ERKC-1x-aux} and $W(t)$ \eqref{Eqn:VV} follows a similar approach with the only difference being in $G_{n,i}^{\rm I}$ and $G^{\rm C}_{n,i}$. 
Specifically, $G^{\rm I}_{n,i}$ incorporates $U(t)$ and its interpolation $\Pi_hU(t)$, while $G^{\rm C}_{n,i}$ is based on $W(t)$ only.
The only difference in the proof comes from the different approximation of the delayed solution: for ERKC-I methods, the error in the delayed solution is estimated with the help of \eqref{Eqn:interpolation-s}, whereas for ERKC-C methods, the error in the delayed solution is estimated directly through $\|W(t)-u(t)\|_V$.
As the proof closely resembles that of ERKC-I methods, it is skipped.

\begin{theorem}\label{Thm:2x-s-global-error}
Let the initial value problem \eqref{Eqn:problem} satisfy Assumptions \ref{Ass:sectorial}-\ref{Ass:V}, and consider for its numerical solution the ERKC-C method \eqref{ERKC-2x}. Let the step size sequence $\{h_j\}_{j=1}^N$ satisfy the condition \eqref{Eqn:step-size}.  If $f:\mathcal{I}^\mu\to X$ \eqref{Eqn:definition-f} is $s$ times differentiable for $\mu=1,2,\ldots,m+1$ with $f^{(s)}\in L^\infty(0,T;X)$, then for $h=\max_{1\leq j \leq N} h_j$ sufficiently small, the numerical solution \eqref{Eqn:VV} is determined uniquely and the error of $W(t)$ satisfies 
\begin{equation*}\label{Eqn:2x-s-V}
\max_{0\leq t \leq T}\|W(t)-u(t)\|_V \leq    Ch^s  \sup_{\sigma\in[0,T]}\|f^{(s)}(\sigma)\|_X   ,
\end{equation*}
where the constant $C$ depends on $T$, but is independent of the step size sequence.
\end{theorem}


\section{Superconvergence}\label{Sec:superconvergence}
It is known that $s$-stage ERKC methods for semilinear parabolic problems (without delay) achieve order $s+1$ at the nodes (see Theorem 5 in \cite{HochbruckOstermann05:323}), provided that the underlying quadrature rule is of order $s+1$, i.e.,
\begin{equation}\label{Eqn:s+1-condition}
\sum_{i=1}^s b_i(0)c_i^s = \frac{1}{s+1},
\end{equation}
and that $f(t)$ has additional smoothness (see Theorem \ref{Thm:ERKCI-s+1}). Condition \eqref{Eqn:s+1-condition} is satisfied, for example, when the collocation parameters are chosen as Radau points with $s\geq 2$ or as Gauss points. In this section, we show that both the ERKC-I and ERKC-C methods with underlying quadrature rule of order $s+1$ actually achieve order $s+1$ over the entire time interval $[0,T]$, rather than only at the nodes.

We shall first extend this superconvergence result to ERKC-I methods for semilinear parabolic problems with time-dependent delay \eqref{Eqn:problem}. 
The following improved error estimate of $\Pi_h U(t)$ in each interval $I_{k+1}$ is required (cf.~\eqref{Eqn:interpolation-s}),
\begin{align}
 \max_{t\in I_{k+1}}\| \Pi_hU (t)-u(t)\|_V  
&  \leq  \max_{t\in I_{k+1}}\|\Pi_hU(t) - \Pi_hu(t) \|_V  +\max_{t\in I_{k+1}}\| \Pi_hu(t) - u(t)\|_V \notag \\
&  \leq C \max_{\sigma \in S_{k+1}} \|U(\sigma)-u(\sigma)\|_V+Ch^{s+1}_{k+1} \sup_{\sigma\in (t_{k},t_{k+1})}\|u^{(s+1)}(\sigma)\|_V, \label{Eqn:interpolation-s+1}
\end{align}
where the constant $C$ is independent of $I_{k+1}$ and the step size sequence.
This is satisfied provided that $u^{(s+1)}\in L^\infty( 0,T;V)$ and $c_1\neq 0$ or $c_s\neq 1$ holds.

The superconvergence result for the ERKC-I method of order $s+1$ is presented in the following theorem.

\begin{theorem}\label{Thm:ERKCI-s+1}
Let the initial value problem \eqref{Eqn:problem} satisfy Assumptions \ref{Ass:sectorial}-\ref{Ass:V}, and consider for its numerical solution the ERKC-I method \eqref{ERKC-1x} whose underlying quadrature rule is of order $s+1$, with collocation parameters $c_1\neq 0$ or $c_s\neq 1$. Let the step size sequence $\{h_j\}_{j=1}^N$ satisfy the condition \eqref{Eqn:step-size}. If $u:\mathcal{I}^\mu\to V$ and $f:\mathcal{I}^\mu\to V$ \eqref{Eqn:definition-f} are both $s+1$ times differentiable for $\mu=1,2,\ldots,m+1$ with $u^{(s+1)}\in L^\infty(0,T;V)$, $f^{(s)}\in L^\infty(0,T;V)$, $f^{(s+1)}\in L^1(0,T;V)$, then for $h=\max_{1\leq j \leq N} h_j$ sufficiently small the error of $\Pi_hU(t)$ satisfies 
$$
\begin{aligned}
\max_{t\in[0,T]} \|\Pi_h U(t)-u(t)\|_{V} & \leq  CC_{\mathrm{S}} h^{s+1} \bigg(\int_0^{T}\|f^{(s+1)}(\sigma)\|_V \mathrm{~d} \sigma+\sup_{\sigma \in [0,T]} \|f^{(s)}(\sigma)\|_V   \bigg)\\
&\quad+ Ch^{s+1} \sup_{\sigma\in[0,T]} \|u^{(s+1)}(\sigma)\|_V.
\end{aligned}
$$
In general, the size of $C_{\mathrm{S}}$ depends on the chosen step size sequence. An upper bound is provided in Lemma \ref{Lem:s+1} below. However, when the step sizes are constant or when the operator $A$ and the space $X$ satisfy certain conditions (see Remark \ref{Rem:Cs} below), $C_{\mathrm{S}}$ is independent of the step size sequence. On the other hand, the constant $C$ depends on $T$, but not on the step size sequence.
\end{theorem}
\begin{proof}
The proof still proceeds by induction and is quite similar to the proof of Theorem \ref{Thm:1x-s-global-error}. We just need an improved estimate of the defects.
Writing the defects \eqref{Eqn:defect} at $\theta=1$ as
\begin{equation}\label{Eqn:split-delta}
\delta_{j}(1) =h_{j}^{s+1}\psi_{s+1}(-h_{j} A)f^{(s)}(t_{j-1}^+) + \widetilde{\delta}_{j}(1),
\end{equation}
with
$$
\psi_{s+1}(-h_{j} A) = \varphi_{s+1}(-h_{j}A) - \frac{1}{s!}\sum_{i=1}^s b_i(-h_{j}A)c_i^s,
$$
we obtain
\begin{equation}\label{Eqn:hat-delta}
\Bigg\|\sum_{j=1}^{n}    \mathrm{e}^{-(t_n-t_j) A}\widetilde{\delta}_j(1) \Bigg\|_V \leq C\sum_{j=1}^{n}  \|  \widetilde{\delta}_j (1) \|_V \leq Ch^{s+1} \int_{0}^{t_n} \|f^{(s+1)}(\sigma)\|_V~\mathrm{d}\sigma.
\end{equation}
In the case of uniform time stepping (i.e., $h_j =h$ for all $j$), from the proof of \cite[Theorem 2]{HochbruckOstermann05:323}, we have
\begin{align}
&\Bigg\| \sum_{j=1}^{n}    h_{j}^{s+1}\mathrm{e}^{-(t_n-t_j) A}{\psi}_{s+1}(-h_{j} A)f^{(s)}(t_{j-1}^+)\Bigg\|_V   \notag \\
&\qquad\leq C_{\mathrm{S}} h^{s+1} \bigg(\int_0^{T}\|f^{(s+1)}(\sigma)\|_V \mathrm{~d} \sigma+\sup_{\sigma \in [0,T]} \|f^{(s)}(\sigma)\|_V   \bigg), \label{Eqn:convolution-tildes+1-Csgood}
\end{align}
where $C_{\mathrm{S}}$ is a constant independent of $n$ and $h$. In the case of variable step sizes satisfying \eqref{Eqn:step-size}, estimate \eqref{Eqn:convolution-tildes+1-Csgood} also holds, provided that additional assumptions on $A$ and the space $X$ are satisfied (see Remark \ref{Rem:Cs}).  In the absence of these assumptions, the following estimate can still be derived from Lemma \ref{Lem:s+1}
\begin{align}
&\Bigg\| \sum_{j=1}^{n}    h_{j}^{s+1}\mathrm{e}^{-(t_n-t_j) A}{\psi}_{s+1}(-h_{j} A)f^{(s)}(t_{j-1}^+)\Bigg\|_V   \leq CC_{\mathrm{S}} h^{s+1}  \sup_{\sigma \in [0,t_{n}] } \|f^{(s)}(\sigma)\|_V    . \label{Eqn:convolution-tildes+1}
\end{align}
Combining \eqref{Eqn:hat-delta} with \eqref{Eqn:convolution-tildes+1-Csgood}-\eqref{Eqn:convolution-tildes+1}, we obtain
$$
\Bigg\|\sum_{j=1}^{n}    \mathrm{e}^{-(t_n-t_j) A} {\delta}_j (1) \Bigg\|_V \leq
CC_{\mathrm{S}} h^{s+1} \bigg(\int_0^{t_n}\|f^{(s+1)}(\sigma)\|_V \mathrm{~d} \sigma+\sup_{\sigma \in [0,t_{n}] } \|f^{(s)}(\sigma)\|_V   \bigg),
$$
which concludes the proof.
\end{proof}

The following lemma was used in the above proof.

\begin{lemma}\label{Lem:s+1}
Let the step size sequence $\{h_j\}_{j=1}^N$ satisfy the condition \eqref{Eqn:step-size} and let $\psi(z)$ be a linear combination of the entire functions  $\varphi_k(z)$ satisfying $\psi(0)=0$. Then, given a sequence $\{w_j\}_{j=0}^{N-1}$ in $V$, the following bound holds for $n=1,2,\ldots,N$,
$$
\Bigg\| \sum_{j=1}^{n}    h_{j}^{s+1}\mathrm{e}^{-(t_n-t_j) A}{\psi}(-h_{j} A)w_{j-1}\Bigg\|_V \leq CC_{\mathrm{S}} h^{s+1} \max_{0\leq j \leq n-1} \|w_j\|_V,
$$
where 
\begin{equation}\label{Eqn:CS}
1\leq C_{\mathrm{S}} \leq \varrho \ln \frac{T}{\min_{1\leq j \leq N} h_j}+2 
\end{equation}
and the constant $C$ is independent of $n$ and the step size sequence.  
\end{lemma}
\begin{proof}
Since $\psi(0)=0$ there exists a bounded operator $\widetilde{\psi}(-h_{n+1}A)$ such that
\begin{equation}\label{Eqn:psi-s+1}
\psi (-h_{n+1}A) = \psi (-h_{n+1}A) -\psi (0) =h_{n+1}A\cdot \widetilde{\psi}(-h_{n+1}A).
\end{equation}
Let
\begin{align}
S(n)& =  \sum_{j=1}^{n}    h_{j}^{s+1}\mathrm{e}^{-(t_n-t_j) A}{\psi} (-h_{j} A)w_{j-1}  . \label{Eqn:Sn}
\end{align}
For $n=1$ and 2, this term is estimated directly as follows
\begin{equation}\label{Eqn:case-(1)-1} 
\|S(n)\|_V \leq C\sum_{j=1}^n h_j^{s+1}\cdot \max_{0\leq j \leq 1} \|w_j\|_V.
\end{equation}
For $n\geq 3$, we replace $\psi_{s+1}$ in \eqref{Eqn:Sn} with \eqref{Eqn:psi-s+1} for $j=1,2,\ldots,n-2$, and get
\begin{equation}\label{Eqn:case-(1)-2} 
\|S(n)\|_V \leq C\Bigg(h^{s+1}\sum_{j=1}^{n-2} h_j(t_n-t_j)^{-1} +h_{n-1}^{s+1}+h^{s+1}_n\Bigg)\max_{0\leq j \leq n-1} \|w_j\|_V.
\end{equation}
Using \eqref{Eqn:step-size}, we derive
\begin{align}
h^{s+1}\sum_{j=1}^{n-2} h_j(t_n-t_j)^{-1}  &\leq \varrho h^{s+1} \sum_{j=1}^{n-2} h_{j+1}(t_n-t_j)^{-1}   \notag \\
& \leq \varrho h^{s+1} \int_{t_1}^{t_{n-1}} (t_n-t)^{-1}~\mathrm{d}t  \notag \\
& = \varrho h^{s+1}   \ln \frac{t_n-t_1}{h_n}  .   \label{Eqn:case-(1)-3} 
\end{align}
A combination of \eqref{Eqn:case-(1)-2}-\eqref{Eqn:case-(1)-3} leads to, for $n\geq 3$,
\begin{equation*}
\|S(n)\|_V \leq CC_{n}h^{s+1} \cdot \max_{0\leq j \leq n-1} \|w_j\|_V,\quad
\mbox{with}~C_{n} = \varrho \ln \frac{t_n-t_1}{h_n}+2.
\end{equation*}
Defining $C_{\mathrm{S}} = \sup_{t_n\leq T} C_n$, it follows that, for $n=1,2,\ldots,N$,
\begin{equation*}
\|S(n)\|_V \leq CC_{\mathrm{S}}h^{s+1} \cdot \max_{0\leq j \leq n-1} \|w_j\|_V,\quad
\mbox{with}~ 1\leq C_{\mathrm{S}}\leq \varrho \ln \frac{T}{\min_{1\leq j \leq N}h_j}+2,
\end{equation*}
which completes the proof.
\end{proof}
\begin{remark}\label{Rem:Cs}
We now discuss the dependence of $C_{\mathrm{S}}$ on the step size sequence in the case of variable time steps satisfying \eqref{Eqn:step-size}. If $X$ is a separable Hilbert space and $A$ is linear, self-adjoint, positive definite with a compact inverse (defined in ${D}(A)\subset X$), by using spectral decomposition and the technique developed in \cite[Theorem 2]{HochbruckOstermann05:323}, one can prove that \eqref{Eqn:convolution-tildes+1-Csgood} is satisfied with $C_{\mathrm{S}}$ independent of $n$ and the step size sequence.

In the worst case, when the step size sequence satisfies \eqref{Eqn:step-size}, the number $C_{\mathrm{S}}$ in estimate \eqref{Eqn:convolution-tildes+1} follows from Lemma~\ref{Lem:s+1}.
In many situations, the number $C_{\mathrm{S}}$ in \eqref{Eqn:CS} can be controlled independently of $\min_{1\leq j \leq N} h_j$. For example, consider a constrained mesh $I_h=\{nh:n=0,1,\ldots,N\} \cup  \{\xi_1,\xi_2,\ldots,\xi_{m}\}$, which is generated by including the potential discontinuity points into the uniform partition. Similar to the derivation of \eqref{Eqn:CS}, it follows that
$\|S(N)\|_V   \leq Ch^{s+1}(\ln N +1)   \max_{0\leq j \leq N-1} \|w_j\|_V$. 
\end{remark}

The superconvergence result for the ERKC-C method of order $s+1$ is presented in the following theorem.

\begin{theorem}\label{Thm:ERKCC-s+1}
Let the initial value problem \eqref{Eqn:problem} satisfy Assumptions \ref{Ass:sectorial}-\ref{Ass:V}, and consider for its numerical solution the ERKC-C method \eqref{ERKC-2x} whose underlying quadrature rule is of order $s+1$. Let the step size sequence $\{h_j\}_{j=1}^N$ satisfy the condition \eqref{Eqn:step-size}. If $f:\mathcal{I}^\mu\to V$ is $s+1$ times differentiable for $\mu=1,2,\ldots,m+1$ with $f^{(s)}\in L^\infty(0,T;V)$, $f^{(s+1)}\in L^1(0,T;V)$, then for $h=\max_{1\leq j \leq N} h_j$ sufficiently small the error of $W(t)$ satisfies 
$$
\begin{aligned}
\max_{t\in[0,T]} \|W(t)-u(t)\|_{V} & \leq  CC_{\mathrm{S}} h^{s+1} \bigg(\int_0^{T}\|f^{(s+1)}(\sigma)\|_V \mathrm{~d} \sigma+\sup_{\sigma \in [0,T]} \|f^{(s)}(\sigma)\|_V   \bigg).
\end{aligned}
$$
In general, the size of $C_{\mathrm{S}}$ depends on the chosen step size sequence. An upper bound is provided in Lemma \ref{Lem:s+1}. However, when the step sizes are constant or when the operator $A$ and the space $X$ satisfy certain conditions (see Remark \ref{Rem:Cs}), $C_{\mathrm{S}}$ is independent of the step size sequence. On the other hand, the constant $C$ depends on $T$, but not on the step size sequence.
\end{theorem}

Theorems \ref{Thm:ERKCI-s+1} and \ref{Thm:ERKCC-s+1} give superconvergence results of order $s+1$. From Theorem 6 in \cite{HochbruckOstermann05:323}, the convergence order of ERKC methods with constant step size further achieves $s+1+\beta$ at the nodes within $[0,\xi_1]$, provided that the underlying quadrature rule is of order $s+2$ and some slightly stronger assumptions regarding spatial regularity hold.
The extension of this result to ERKC-I and ERKC-C methods over $[0,T]$ is not feasible under the present framework. 
To address this limitation, we propose the following modified ERKC-I method by adjusting the interpolation approach, 
\begin{equation}\label{mERKC-1x}
\begin{aligned}
U_{n+1}
& = \mathrm{e}^{-h_{n+1}A}U_n + h_{n+1} \sum_{i=1}^s b_i(-h_{n+1}A) G^{\rm I}_{n,i}, \\
U_{n,i} &
= \mathrm{e}^{-c_i h_{n+1}A}U_n + h_{n+1} \sum_{j=1}^s a_{ij}(- h_{n+1}A) G^{\rm I}_{n,j},\quad 1 \leq i \leq s, \\
G^{\rm I}_{n,j}  & =g(t_{n,j},U_{n,j},(\Pi_h^{*}U)^\tau(t_{n,j})),\quad 1\leq j \leq s.
\end{aligned}
\end{equation}
For $t\in [t_n,t_{n+1}]$, we have that $\hat{t}=t-\tau(t)$ lies either in $[-\tau_0,0]$ or in $\mathcal{I}^\mu$ for some $0\leq\mu\leq m$. We set $\Pi^*_h U(\hat{t})=\phi(\hat{t})$ for $\hat{t}\in[-\tau_0,0]$. 
For $\hat{t}\in I_{k+1}\subseteq \mathcal{I}^\mu$, the value of $\Pi^*_h U(\hat{t})$ is given by a polynomial based on the interpolation at mesh points as follows:
$$
\Pi_h^* U(\hat{t}) = P_{k+1}(\hat{t}),\quad \hat{t} \in I_{k+1}\subseteq \mathcal{I}^\mu,
$$
where $P_{k+1}$ is a $(\ell_k+r_k-1)$-degree interpolation polynomial (with $\ell_k,r_k\geq 1$) on $I_{k+1}$ satisfying
\begin{equation}\label{Eqn:modified-PiU}
P_{k+1} (t_{q})=U_{q} \quad \mbox{for}~t_q \in S^{*}_{k+1}=\{t_{k-\ell_k+1},\ldots,t_k,t_{k+1},\ldots,t_{k+r_k}\}\subseteq \mathcal{I}^\mu . 
\end{equation}
In practice, a symmetric (or nearly symmetric, for $s$ being odd) strategy based on $s+2$ interpolation points is used for choosing $\ell_k$ and $r_k$, with the necessary adjustments to ensure that  the interpolation stencil does not cross primary discontinuities.

If $\ell_k+r_k=s+2$ and $u^{(s+2)}\in L^\infty(0,T;V)$, the error of the numerical solution $\Pi_h^* U(t)$ in each interval $I_{k+1}$ can be controlled by 
\begin{align}
\max_{t\in I_{k+1}}\| \Pi_h^*U (t)-u(t)\|_V \notag 
& \leq  \max_{t\in I_{k+1}}\|\Pi_h^*U(t) - \Pi_h^*u(t) \|_V  +\max_{t\in I_{k+1}}\| \Pi_h^*u(t) - u(t)\|_V \notag \\
&  \leq C \max_{\sigma \in S_{k+1}^*} \|U(\sigma)-u(\sigma)\|_V+Ch^{s+2} \sup_{\sigma\in (t_{k-\ell_k+1},t_{k+r_k})}\|u^{(s+2)}(\sigma)\|_V. \label{Eqn:interpolation-s+2}
\end{align}
The number $C$ is independent of $I_{k+1}$ and the step size sequence provided that the following condition holds: there are constants $\kappa_1$ and $\kappa_2$ such that
\begin{equation}\label{Eqn:k1k2}
\kappa_1\leq h_{j+1}/h_j \leq \kappa_2\quad\mbox{for all}~ j.
\end{equation}
In particular, quasi-uniform meshes satisfy this condition.

Moreover, we need the following assumptions on spatial regularity.

\begin{assumption}\label{Ass:s+2}
Let $g:[0, T] \times V\times V \rightarrow X$ be Fr\'echet differentiable with respect to the second variable and let $ \frac{\partial g}{\partial v} (t,v,w)$ be locally Lipschitz-continuous in a strip along the exact solution. Further, let $\beta\in [ 0,1]$ be such that
\begin{align}
&\Bigg\| A^\beta \frac{\partial g}{\partial u} (t,u(t),u^\tau(t)  )w\Bigg\|_X
\leq C\|A^\beta w\|_V, \quad \forall w\in {D}(A^{\alpha+\beta}),\label{Eqn:ass-s+2-1}
\end{align}
uniformly for $0\leq t\leq T$.
\end{assumption}

 From Assumption \ref{Ass:s+2} we infer that there exists a real number $\widetilde{L}(R,T)$ such that 
\begin{equation}\label{Eqn:improved-lip}
\begin{aligned}
   \bigg\|g(t, v,u^\tau(t))-g(t,u(t), u^\tau(t))-\frac{\partial g}{\partial u }(t,u(t), u^\tau(t))(v-u(t))   \bigg\|_X  \leq  \widetilde{L} \|v-u(t) \|_V^2   ,
\end{aligned}
\end{equation}
for all $t \in[0, T]$ and all $v \in V$ with $\|v-u(t)\|_V \leq R$.  

The superconvergence result for the ERKC-I method of order $s+1+\beta$ is presented in the following theorem.

\begin{theorem}\label{Thm:ERKC-I-s+2}
Let the initial value problem \eqref{Eqn:problem} satisfy Assumptions \ref{Ass:sectorial}-\ref{Ass:V} and \ref{Ass:s+2}, and consider for its numerical solution the modified ERKC-I method \eqref{mERKC-1x} whose underlying quadrature rule is of order $s+2$. Let the interpolation approach \eqref{Eqn:modified-PiU} satisfy $\ell_k+r_k=s+2$ for $k=0,1,\ldots,N-1$ and let the step size sequence $\{h_j\}_{j=1}^N$ satisfy the condition \eqref{Eqn:k1k2}.
If $u:\mathcal{I}^\mu\to V$ and $f:\mathcal{I}^\mu\to V$ \eqref{Eqn:definition-f} are both $s+2$ times differentiable  for $\mu=1,2,\ldots,m+1$ with $u^{(s+2)}\in L^\infty( 0,T;V)$, $A^\beta f^{(s)}  \in L^\infty(0,T;V)$, $f^{(s+1)}\in L^\infty(0,T;V)$ and $f^{(s+2)}\in L^1(0,T;V)$, then for $h=\max_{1\leq j \leq N} h_j$ sufficiently small the error of $U_n$ satisfies 
$$
\max_{0\leq n \leq N} \|U_n-u(t_n)\|_V \leq CC_{\mathrm{S}}h^{s+1+\beta},
$$
where $C_{\mathrm{S}}$ is given by Lemma \ref{Lem:s+1}. The constant $C$ depends on $T$, but is independent of the step size sequence.
\end{theorem}
\begin{proof}
We start with rewriting the error recursion \eqref{Eqn:main-error-1x} as
\begin{equation}\label{Eqn:rewrite-recursion}
\begin{aligned}
e(t_{n+1}) & = \mathrm{e}^{-h_{n+1}A}e(t_n) + h_{n+1}\sum_{i=1}^s b_i(-h_{n+1}A)J_ne(t_{n,i}) \\
&\quad + h_{n+1}\sum_{i=1}^sb_i(-h_{n+1}A)\Big( g(t_{n,i},U_{n,i},(\Pi^*_hU)^{\tau}_{n,i}) - f(t_{n,i}) - J_n e(t_{n,i}) \Big)  - \delta_{n+1}(1),   
\end{aligned}
\end{equation}
with
$$
J_n = \frac{\partial g}{\partial u}(t_n,u(t_n),u^\tau(t_n)).
$$
Solving the above error recursion for $t_n\in (0,\xi_1)$ gives
\begin{align}
e(t_{n}) & =  \mathrm{e}^{-(t_n-t_0)A}  e(t_0)+  \sum_{j=1}^{n}\mathrm{e}^{-(t_n-t_j)A}\bigg( \sum_{i=1}^s h_{j} b_i(-h_{j}A)J_{j-1}e(t_{j-1,i}) \notag \\
&\quad + h_{j}\sum_{i=1}^sb_i(-h_{j}A)\Big( g(t_{j-1,i},U_{j-1,i},(\Pi^*_hU)^{\tau}_{j-1,i})   - f(t_{j-1,i}) - J_{j-1} e(t_{j-1,i}) \Big) - \delta_{j}(1) \bigg). \label{Eqn:solving-error-recursion}
\end{align}
The defects $\delta_{n+1}(1)$ are represented as
$$
\delta_{n+1}(1) = h_{n+1}^{s+1}\psi_{s+1}(-h_{n+1}A)f^{(s)}(t_n^+) +  h_{n+1}^{s+2}\psi_{s+2}(-h_{n+1}A)f^{(s+1)}(t_n^+) +\hat{\delta}_{n+1}(1)
$$
with $\psi_{s+1}(0)=0$ and $\psi_{s+2}(0)=0$ since the underlying Runge--Kutta method has order $s+2$. Consequently, the second and the third term can be bounded as in Theorem \ref{Thm:ERKCI-s+1}. The first term with $\psi_{s+1}$ is rewritten as 
$$
\begin{aligned}
h_{n+1}^{s+1}\psi_{s+1}(-h_{n+1}A)f^{(s)}(t_n^+) &=h_{n+1}^{s+1}\big(\psi_{s+1}(-h_{n+1}A)-\psi_{s+1}(0) \big)f^{(s)}(t_n^+)\\
& =  h_{n+1}^{s+1+\beta}\psi_{s+1}^{(1)}(-h_{n+1}A)(h_{n+1}A)^{1-\beta}\cdot A^\beta f^{(s)}(t_n^+)
\end{aligned}
$$
with $\psi_{s+1}^{(1)}(0)=0$, which can also be bounded in the desired way by using Lemma~\ref{Lem:s+1}. The convolution of $\mathrm{e}^{-(t_n-t_j)A}$ with the defects $\delta_{j}(1)$ is thus bounded by
\begin{equation}\label{Eqn:s+2-defects-node}
\Bigg\|\sum_{j=1}^{n}    \mathrm{e}^{-(t_n-t_j) A} {\delta}_j (1) \Bigg\|_V \leq CC_{\mathrm{S}} h^{s+1+\beta}.
\end{equation}

The defects $\delta_{n+1}(c_i)$ of the inner stages can be rewritten as
\begin{equation*}\label{Eqn:rewritten-internal-defects}
\delta_{n+1}(c_i) = h_{n+1}^{s+1}\psi_{i,s+1}(-h_{n+1}A)f^{(s)}(t_n^+)+\widetilde{\delta}_{n+1}(c_i),
\end{equation*}
with
$$
\psi_{i,s+1}(-h_{n+1}A)=\varphi_{s+1}(-c_ih_{n+1}A)c_i^{s+1}-\frac{1}{s!}\sum_{j=1}^s a_{ij}(-h_{n+1}A)c_j^s.
$$
Following the idea presented in the proof of Theorem 6 in \cite{HochbruckOstermann05:323}, with the help of the relation
$$
\sum_{i=1}^s b_i(0)\psi_{i,s+1}(0)=0
$$
(since the underlying Runge--Kutta method has order $s+2$), the convolution associated with the defects of the inner stages $\delta_j(c_i)$ is bounded by
\begin{equation}\label{Eqn:s+2-defects-innerstages}
\Bigg\|\sum_{j=1}^{n}    \mathrm{e}^{-(t_n-t_j) A} h_{j}\sum_{i=1}^s b_i(-h_{j}A)J_{j-1} \delta_{j}(c_i) \Bigg\|_V \leq C h^{s+1+\beta}.
\end{equation}
Replacing $e(t_{j-1,i})$ in \eqref{Eqn:solving-error-recursion} with \eqref{Eqn:main-error-1x}, the dominant error term comes from 
\begin{equation}\label{Eqn:domaint-error-term}
h_j \sum_{i=1}^s b_i(-h_j A) J_{j-1} \cdot h_j \sum_{m=1}^s a_{i m}(-h_j A) J_{j-1} \delta_{j}(c_m).
\end{equation}
We can derive the desired bound for the convolution of $\mathrm{e}^{-(t_n-t_j)A}$ with \eqref{Eqn:domaint-error-term} in the same manner as  presented in \cite{HochbruckOstermann05:323}.
Using \eqref{Eqn:s+2-defects-node}, \eqref{Eqn:s+2-defects-innerstages} and \eqref{Eqn:improved-lip}, we get by applying the discrete Gronwall lemma (Lemma \ref{Lem:gronwall}) the following result: for any mesh point $t_n$ belonging to $[0,\xi_{1}]$, it holds
\begin{equation*}
\|e(t_n)\|_{V}  \leq CC_{\mathrm{S}} h^{s+1+\beta}.
\end{equation*}
Combining this with \eqref{Eqn:interpolation-s+2}, it holds
$$
\max_{t\in[0,\xi_{1}]} \|\Pi^*_hU(t) - u(t)\|_V \leq CC_{\mathrm{S}}h^{s+1+\beta}.
$$

The proof again proceeds by induction. We will show that if
\begin{equation}
\begin{aligned}
\max_{t\in[0,\xi_l]} \|\Pi^*_hU(t)-u(t)\|_V & \leq CC_{\mathrm{S}}h^{s+1+\beta}
\label{Eqn:s+2-assume-0k} 
\end{aligned}
\end{equation}
holds for $l= \mu$, then it also holds for $l=\mu+1$. 
By \eqref{Eqn:improved-lip}, we have
\begin{align}
& \big\| g\big(t_{j-1,i},U_{j-1,i},(\Pi^*_hU)^{\tau}_{j-1,i}\big) - f(t_{j-1,i}) - J_{j-1} e(t_{j-1,i}) \big\|_X \notag \\
&\quad \leq \Big\| g\big(t_{j-1,i},U_{j-1,i},u^{\tau}(t_{j-1,i})\big) - f(t_{j-1,i}) - \frac{\partial g}{\partial u}\big(t_{j-1,i},u(t_{j-1,i}),u^\tau(t_{j-1,i})\big) e(t_{j-1,i}) \Big\|_X \notag \\
&\qquad + \big\|g\big(t_{j-1,i},U_{j-1,i},(\Pi^*_hU)^{\tau}_{j-1,i}\big)    - g\big(t_{j-1,i},U_{j-1,i},u^{\tau}(t_{j-1,i})\big)\big\|_X  \notag \\
&\qquad+ \Big\|\frac{\partial g}{\partial u}\big(t_{j-1,i},u(t_{j-1,i}),u^\tau(t_{j-1,i}) \big)e(t_{j-1,i})- J_{j-1} e(t_{j-1,i}) \Big\|_X
\notag \\
&\quad \leq  \widetilde{L} \|e(t_{j-1,i}) \|_V^2 
 + L\|(\Pi_h^*U)^{\tau}_{j-1,i}-u^\tau(t_{j-1,i})\|_V+Ch\|e(t_{j-1,i})\|_V.\label{Eqn:s+2-nonlinearity}
\end{align}
We solve the error recursion \eqref{Eqn:rewrite-recursion} for $t_n\in(\xi_{\mu},\xi_{\mu+1})$ and consider its estimate in $V$. By using \eqref{Eqn:s+2-defects-node}, \eqref{Eqn:s+2-defects-innerstages}, \eqref{Eqn:s+2-assume-0k}  and \eqref{Eqn:s+2-nonlinearity}, we get by applying the discrete Gronwall lemma (Lemma \ref{Lem:gronwall}) the following result: for any mesh point $t_n$ belonging to $[\xi_{\mu},\xi_{\mu+1}]$, it holds
\begin{equation*}
\|e(t_n)\|_{V}  \leq CC_{\mathrm{S}} h^{s+1+\beta}.
\end{equation*}
Combining this estimate with \eqref{Eqn:interpolation-s+2}, it holds
$$
\max_{t\in[0,\xi_{\mu+1}]} \|\Pi^*_hU(t) - u(t)\|_V \leq CC_{\mathrm{S}}h^{s+1+\beta},
$$
and the proof is completed.
\end{proof}
As mentioned in \cite{HochbruckOstermann05:323}, the restriction $\beta\leq 1$ was made just for simplicity. Higher order of convergence can be shown provided that the solution has higher spatial regularity and the underlying quadrature rule is of higher order as well. 
Particularly, full order of convergence can be obtained for problems under periodic boundary conditions.

\section{Numerical experiments}\label{Sec:experiments}
This section is devoted to illustrate the convergence results obtained in the previous sections. Specifically, we consider the problem \eqref{Eqn:problem} with $A=-\Delta$, defined on the computational domain $\Omega=[0,1]^d$, $d=1,2$.

\begin{example}\label{Exa:1}\rm
We begin by investigating the following one-dimensional parabolic problem with known exact solution
\begin{equation}\label{Eqn:example-1}
\frac{\partial u}{\partial t} (t,x) -\frac{\partial^2 u}{\partial x^2}(t,x) = \frac{1}{1+u(t,x)^2}+\frac{1}{1+u(\frac{t}{2}-\frac{1}{2},x)^2}+\Phi(t,x).
\end{equation}
We consider this problem for $x\in[0,1]$ and $t\in [0,3]$, subject to homogeneous Dirichlet boundary conditions. The source function $\Phi$ is determined by the exact solution of the problem $u(t,x)= \Psi(t) \sin x \sin(1-x)$, where $\Psi(t)$ is given by
$$
\Psi(t) = 
\left\{\begin{aligned}
& \mathrm{e}^{-t}, && t\in[-\tfrac{1}{2},0], \\
& 1+t\mathrm{e}^{2t}, && t\in(0,1], \\
& (1+\mathrm{e}^{2})+3\mathrm{e}^{2}(t-1)+(t-1)^2\mathrm{e}^{3t}, && t \in(1,3]
\end{aligned}\right.
$$
such that $\Psi^\prime(0^{-}) \neq \Psi^\prime(0^{+})$ and $\Psi^{\prime\prime}(1^{-}) \neq \Psi^{\prime\prime}(1^{+})$, which meets the characteristics of delay differential equations, considering that primary discontinuities are common in such equations.  
\end{example}

\begin{figure}[H]
    \centering
    \includegraphics[width=0.49\textwidth]{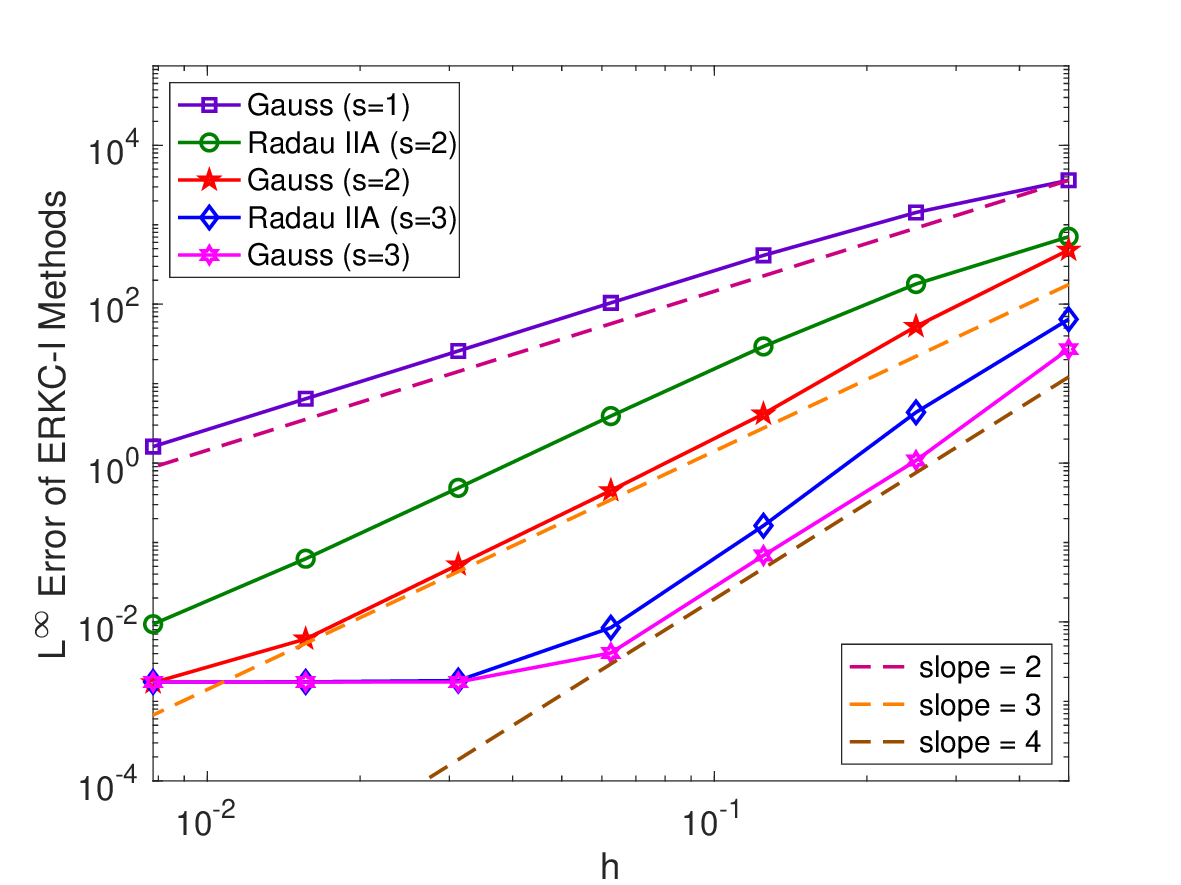}
    \includegraphics[width=0.49\textwidth]{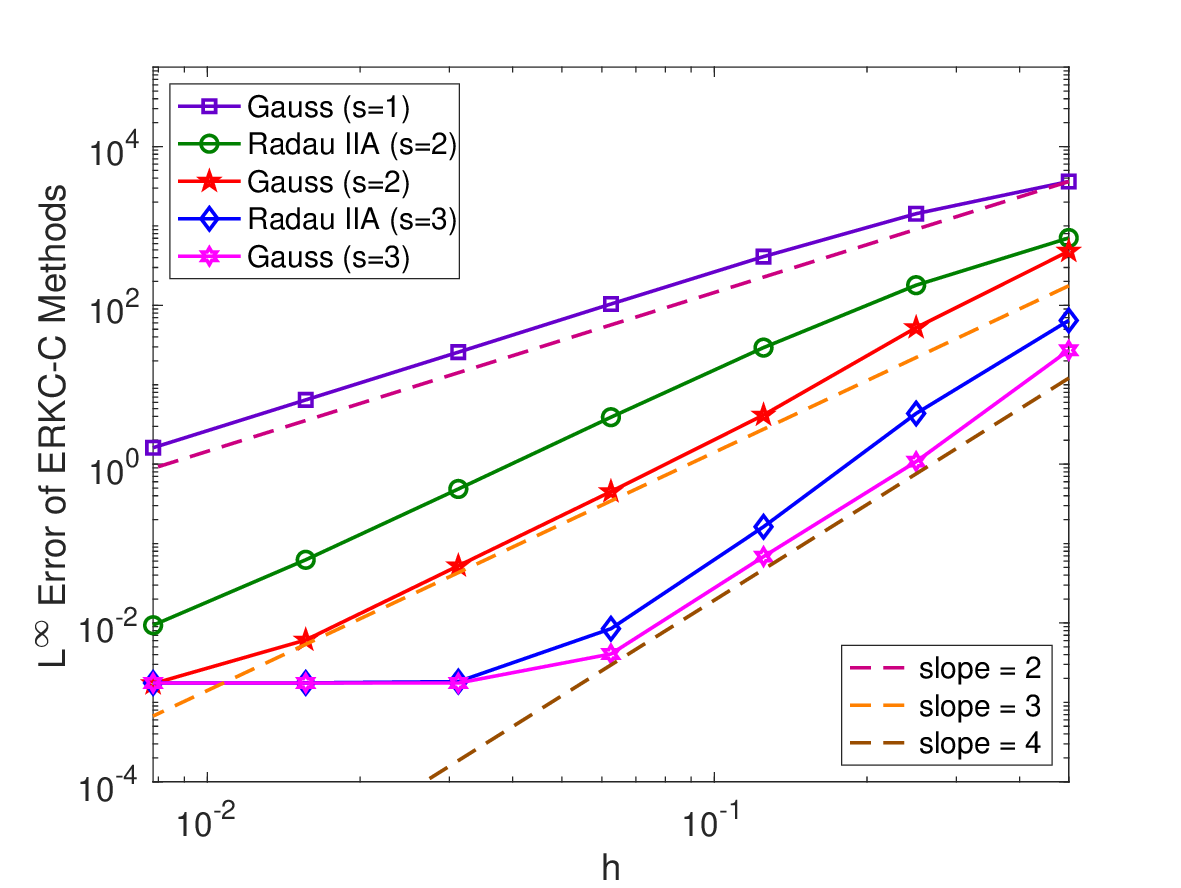}
    \caption{The convergence rates of ERKC-I methods (in the left panel) and ERKC-C methods (in the right panel) for \eqref{Eqn:example-1}. The errors are measured at $T=3$ in the $L^\infty(\Omega)$ norm. In this example, the two considered approaches produced nearly identical errors. Therefore, the two panels are indistinguishable.}
    \label{Fig:source-EndL00} 
\end{figure}

\begin{figure}[H]
    \centering
    \includegraphics[width=0.49\textwidth]{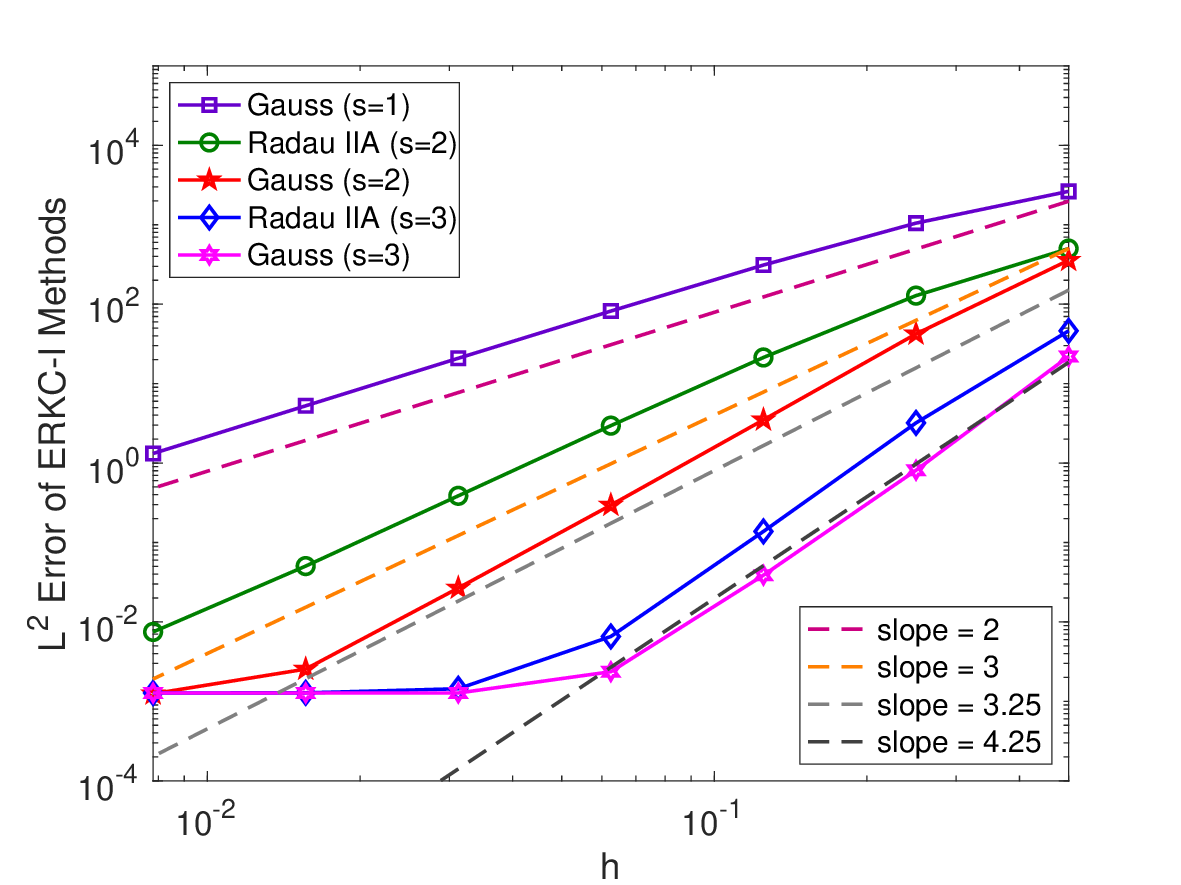} 
    \includegraphics[width=0.49\textwidth]{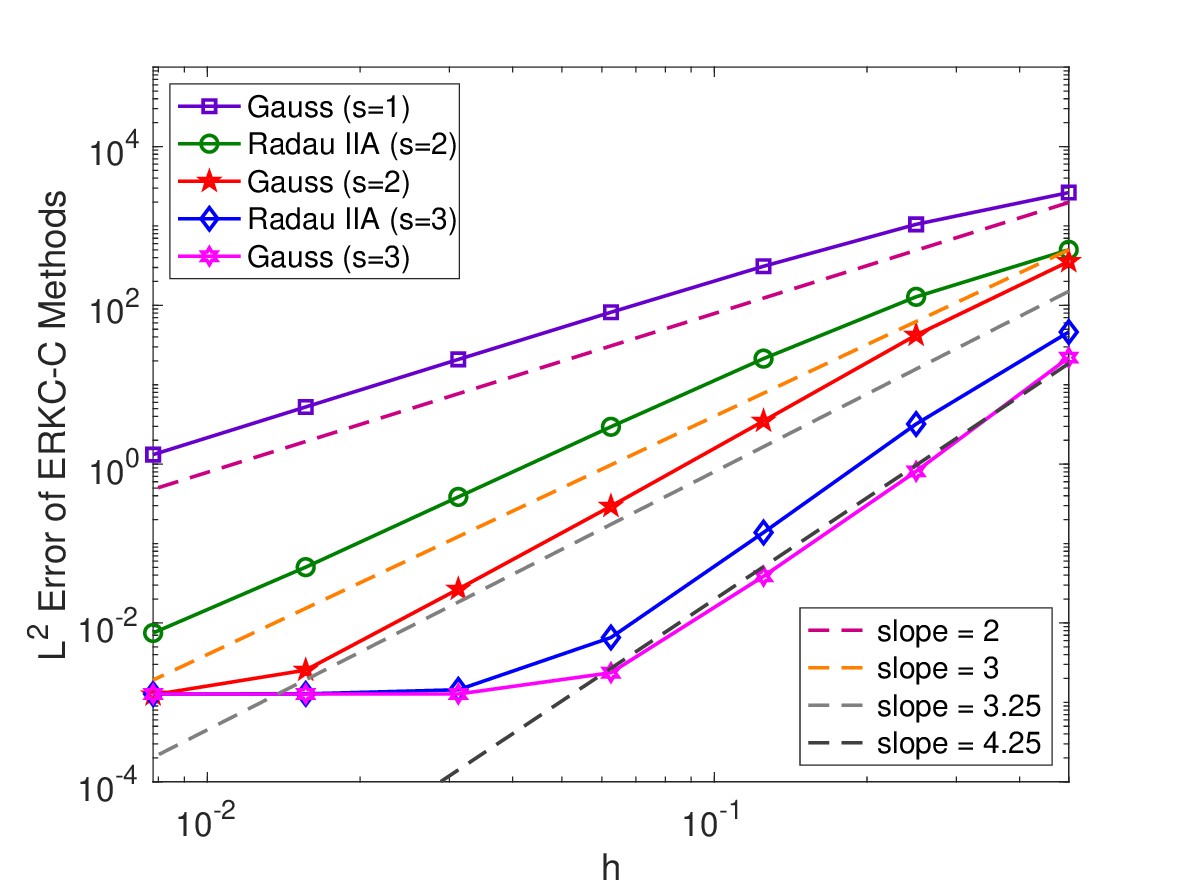}
    \caption{The convergence rates of ERKC-I methods (in the left panel) and ERKC-C methods (in the right panel) for \eqref{Eqn:example-1}. The errors are measured at $T=3$ in the $L^2(\Omega)$ norm. In this example, the two considered approaches produced nearly identical errors. Therefore, the two panels are indistinguishable.}
    \label{Fig:source-EndL2} 
\end{figure}

We apply standard finite differences with $n=1000$ grid points to spatially discretize the problem \eqref{Eqn:example-1}. The resulting products of matrix functions with vectors are computed by the fast Fourier transformation.
The convergence rates of ERKC-I and ERKC-C methods at $T=3$ in the $L^\infty(\Omega)$  and the $L^2(\Omega)$ norms, respectively, are presented in Figures \ref{Fig:source-EndL00} and \ref{Fig:source-EndL2}. The numerical results clearly match our theoretical analysis. 

 In this example, we also observe that both ERKC-I and ERKC-C methods achieve a superconvergence rate of $s+1+\beta$ in the  $L^2(\Omega)$ norm, despite the proof of this superconvergence result being limited to the modified ERKC-I method (with $X=L^2(\Omega)$, the best attainable value of $\beta$ in Theorem \ref{Thm:ERKC-I-s+2} is $\beta=1/4-\varepsilon$ for arbitrary $\varepsilon>0$, see \cite{Fujiwara67:82}). Nevertheless, this superconvergence result generally does not hold for ERKC-I and ERKC-C methods (cf. Example \ref{Exa:mERKC-I}).

\begin{example}\label{Exa:2D}\rm
Consider the two-dimensional semilinear parabolic problem
\begin{equation}\label{Eqn:example-2D}
\frac{\partial u}{\partial t} (t,x,y) -\frac{\partial^2 u}{\partial x^2}(t,x,y)-\frac{\partial^2 u}{\partial y^2}(t,x,y) = \frac{1}{1+u(t,x,y)^2 }+\frac{1}{1+ u(\frac{t}{2}-\frac{1}{2},x,y)^2},
\end{equation}
for $(x,y)\in[0,1]^2$ and $t\in [0,3]$, subject to homogeneous Dirichlet boundary conditions. The initial condition is given by $\phi(t,x,y)=\mathrm{e}^{-t}x(1-x)y(1-y)$ for $t\in [-\frac{1}{2},0]$.
\end{example}

\begin{figure}[H]
    \centering
    \includegraphics[width=0.49\textwidth]{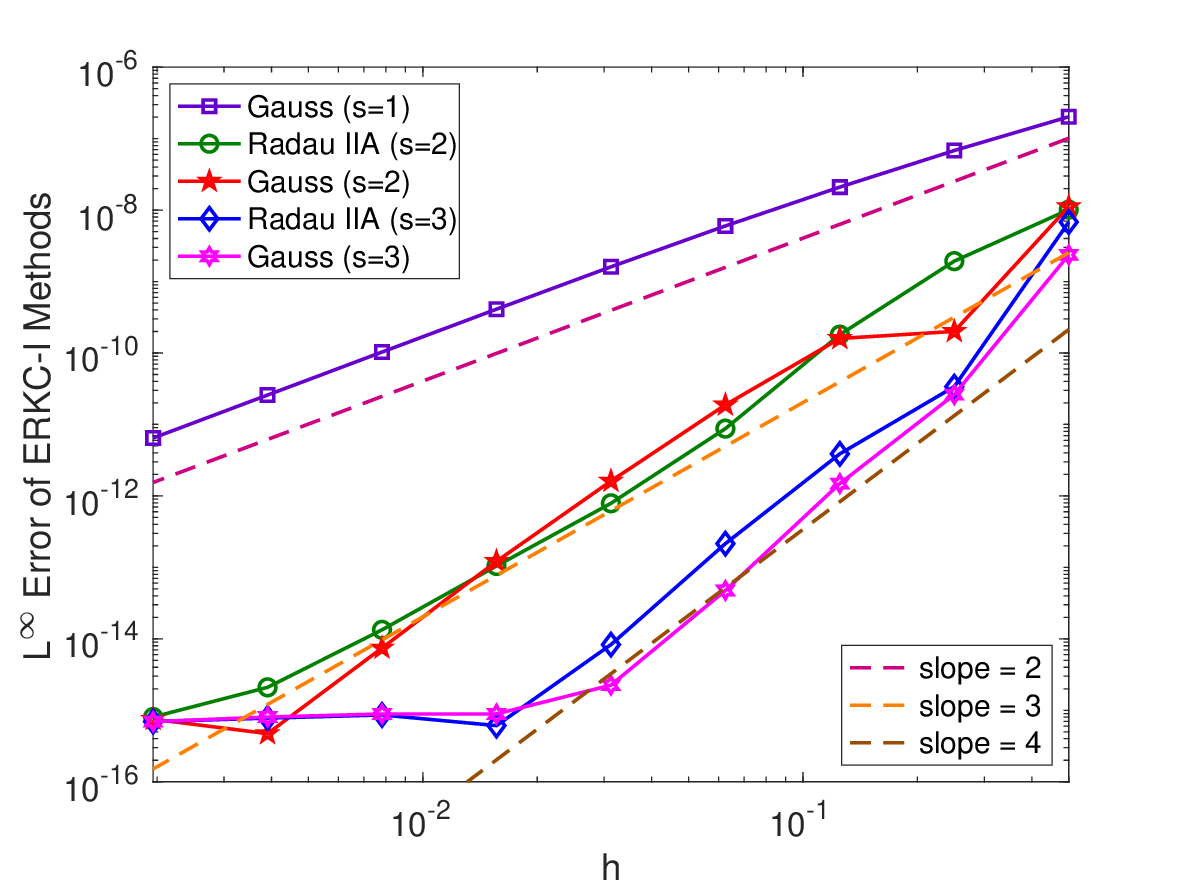}
    \includegraphics[width=0.49\textwidth]{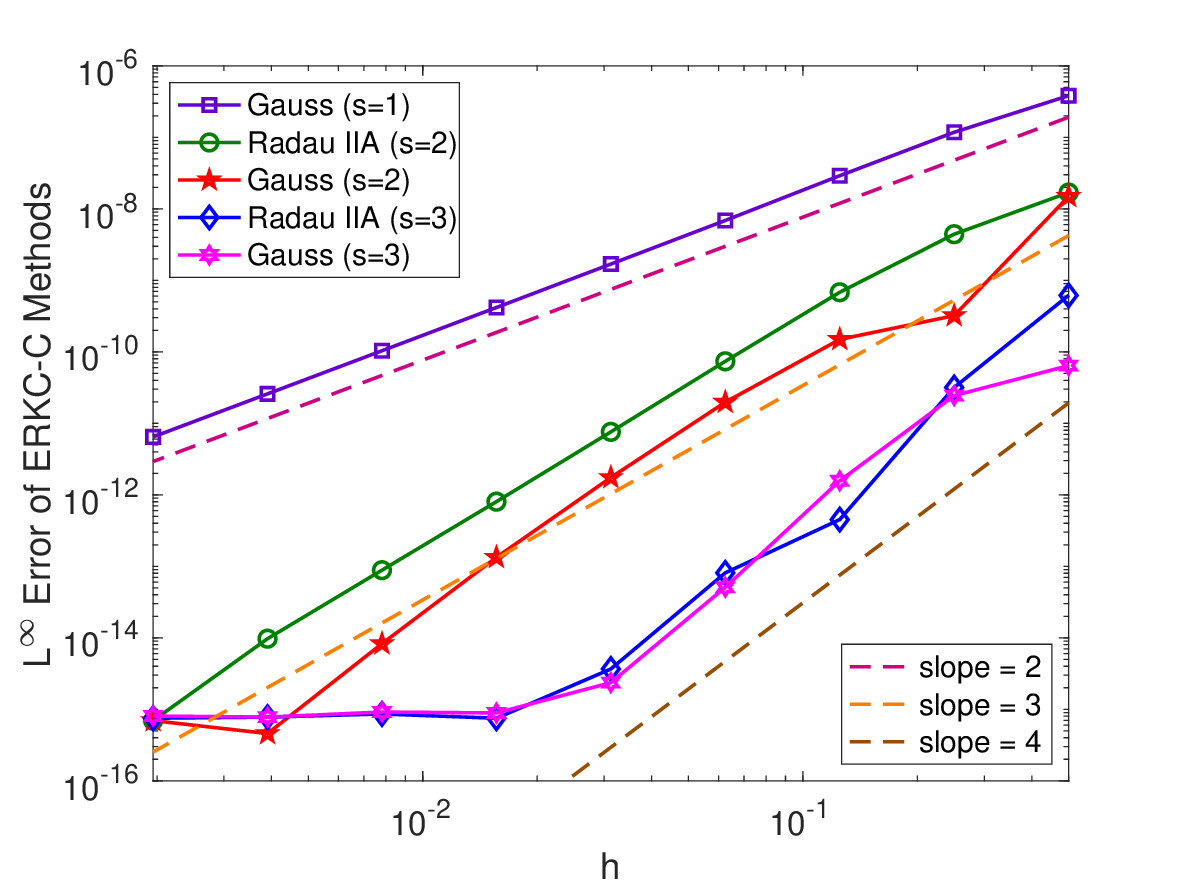}
    \caption{The convergence rates of ERKC-I methods (in the left panel) and ERKC-C methods (in the right panel) for \eqref{Eqn:example-2D}. The errors are measured at $T=3$ in the $L^\infty(\Omega)$ norm.}
    \label{Fig:2Dnosource-EndL00} 
\end{figure}

\begin{figure}[h]
    \centering
    \includegraphics[width=0.49\textwidth]{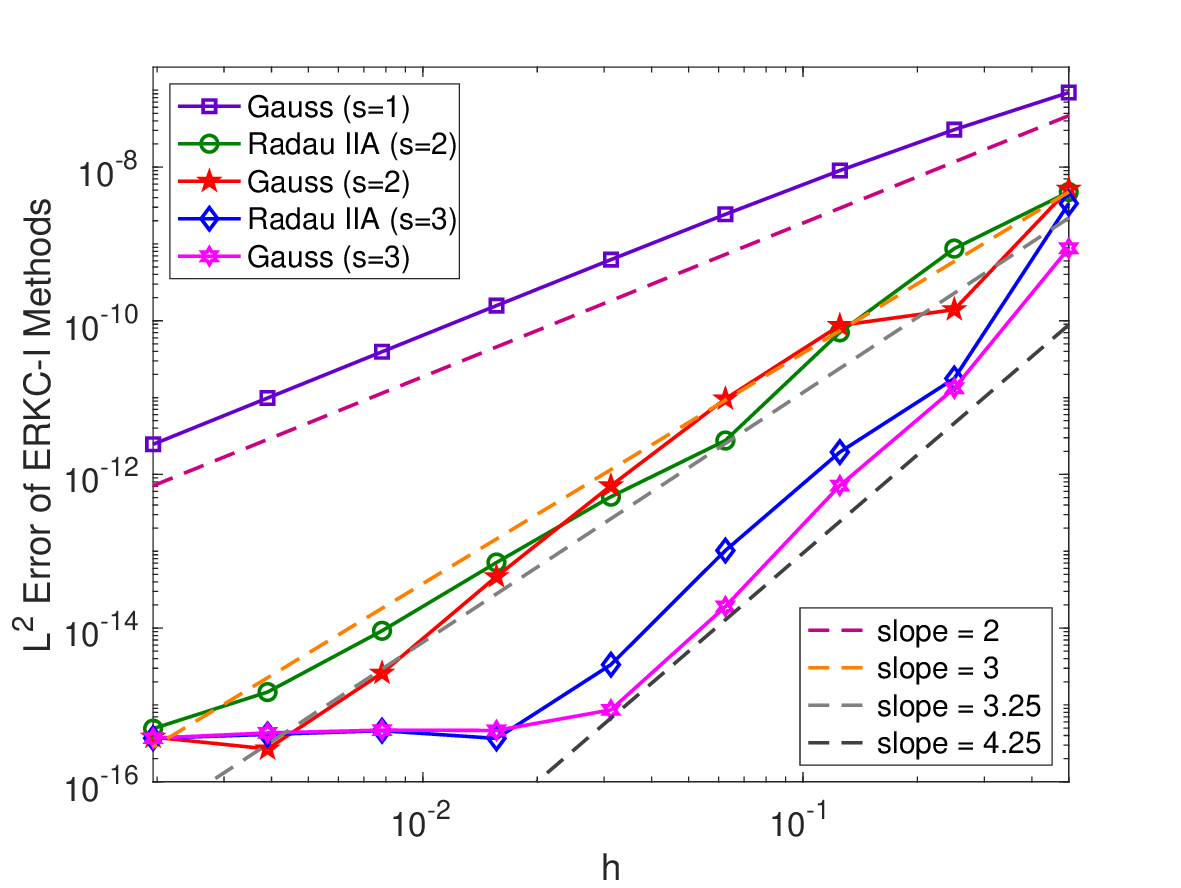}
    \includegraphics[width=0.49\textwidth]{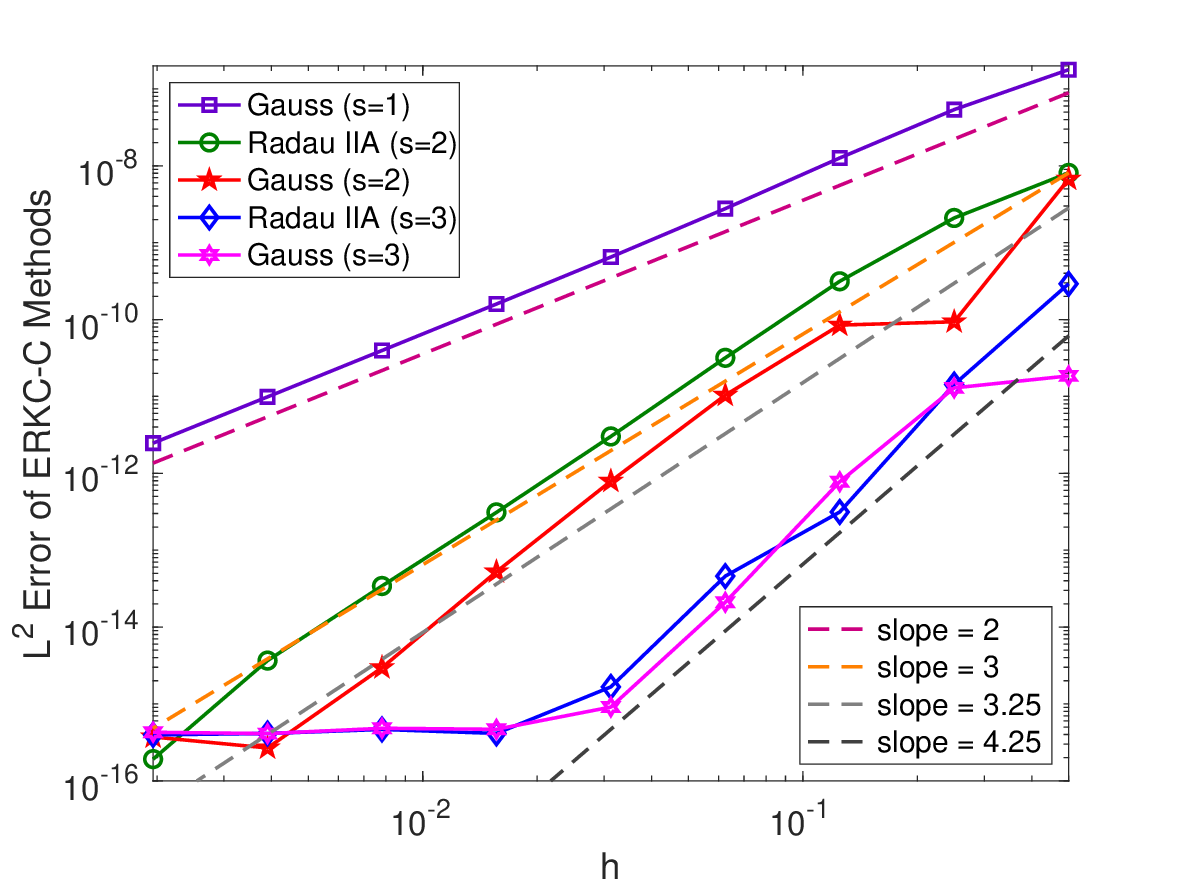}
    \caption{The convergence rates of ERKC-I methods (in the left panel) and ERKC-C methods (in the right panel) for \eqref{Eqn:example-2D}. The errors are measured at $T=3$ in the $L^2(\Omega)$ norm.}
    \label{Fig:2Dnosource-EndL2} 
\end{figure}

We apply standard finite differences with $n=200$ grid points to discretize the problem in each spatial direction. In this example the exact solution is unknown. The reference solution is computed by the ERKC-C method with Gauss collocation points ($s=3$) using the constant step size $h=2^{-11}$. The resulting mesh includes the primary discontinuities.    
The errors of the ERKC-I and ERKC-C methods in the $L^\infty(\Omega)$ and the $L^2(\Omega)$ norm at final time $T=3$ are presented in Figures \ref{Fig:2Dnosource-EndL00} and  \ref{Fig:2Dnosource-EndL2}. The figures confirm the theoretical analysis.  Moreover, a comparison of computational cost of the ERKC methods with Gauss collocation points ($s=3$) is presented in Table~\ref{Tab:comparison}. For step sizes $h=2^{-k}$ with $k=3,4,\ldots,9$, the ERKC-I method consistently reduces the CPU time by approximately 12\% compared to the ERKC-C method, demonstrating its computational efficiency across a range of step sizes.

\begin{table}[h]
\caption{CPU time comparison of ERKC methods with $s=3$ for \eqref{Eqn:example-2D}.}\label{Tab:comparison}
\centering
\renewcommand{\arraystretch}{1.2} 
\setlength{\tabcolsep}{5pt} 
\begin{tabular}{lccccccccc}
\toprule
 & \multicolumn{9}{c}{Step size $h$} \\ 
\cmidrule(lr){2-10} 
Method & $2^{-1}$ & $2^{-2}$ & $2^{-3}$ & $2^{-4}$ & $2^{-5}$ & $2^{-6}$ & $2^{-7}$ & $2^{-8}$ & $2^{-9}$ \\
\midrule
ERKC-I & 10.71 & 15.40 & 27.79 & 52.49 & 92.80 & 179.52 & 351.54 & 662.33 & 1278.92 \\
ERKC-C & 9.90 & 16.82 & 31.66 & 60.21 & 107.64 & 209.64 & 410.12 & 744.17 & 1454.16 \\
\bottomrule
\end{tabular}
\end{table}

\begin{example}\label{Exa:H14}\rm
In this example, we consider the semilinear parabolic problem
\begin{equation}\label{Eqn:H14}
\frac{\partial u}{\partial t} (t,x) -\frac{\partial^2 u}{\partial x^2}(t,x) = u(t,x)(1-u(t,x)) - u\Big(\frac{t}{2}-\frac{1}{2},x\Big)\Big(1-u\Big(\frac{t}{2}-\frac{1}{2},x\Big)\Big)  ,
\end{equation}
for $x\in[0,1]$ and $t\in [0,3]$, subject to homogenous Dirichlet boundary conditions. The initial condition is given by $\phi(t,x)=\mathrm{e}^{t}x(1-x)$ for $t\in [-\frac{1}{2},0]$. The above initial-boundary value problem can be written as an abstract initial value problem with linear operator $A$ defined by $(Au)(x)=-\partial_{xx}u(x)$, and nonlinearity $g$ defined by $g(t,v,w)=v(1-v)-w(1-w)$. The  Banach space $X$ can be chosen as  $X=L^2(\Omega)$ and then $D(A)= H_0^1(\Omega)\cap H^2(\Omega)$. The nonlinearity $g:[0,3]\times V\times V\to X$ is locally Lipschitz-continuous in a strip along the exact solution for $V=H^{\frac{1}{4}}(\Omega)$.
\end{example}

\begin{figure}[H]
    \centering
    \includegraphics[width=0.49\textwidth]{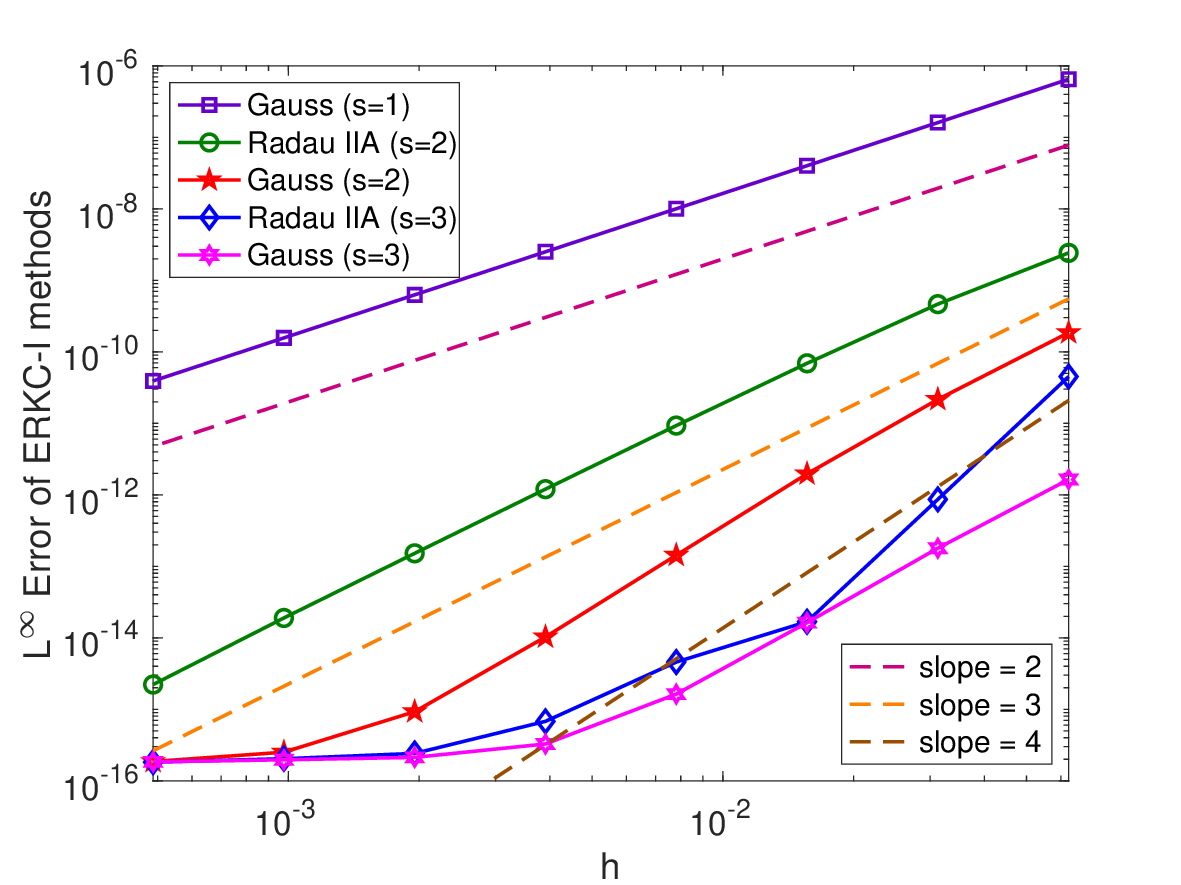}
    \includegraphics[width=0.49\textwidth]{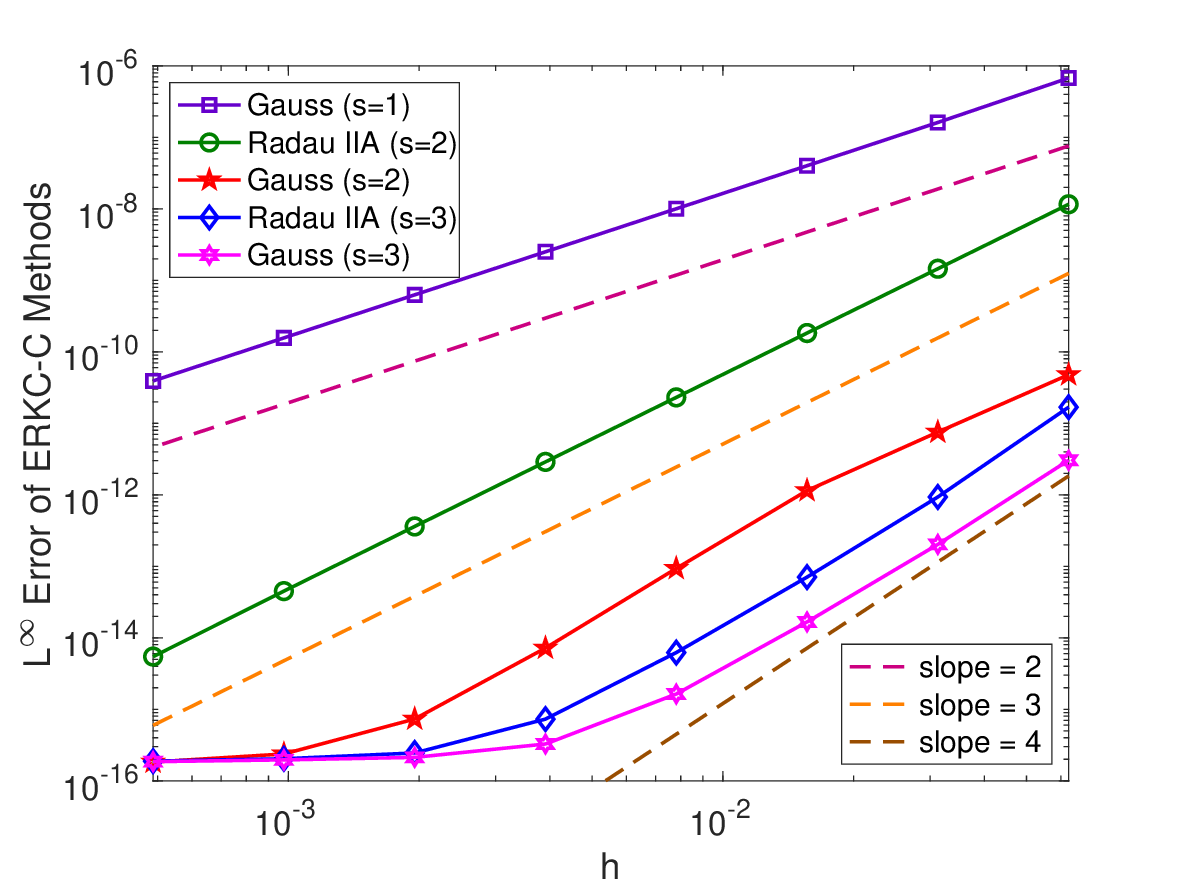}
    \caption{The convergence rates of ERKC-I methods (in the left panel) and ERKC-C methods (in the right panel) for \eqref{Eqn:H14}. The errors are measured at $T=3$ in the $L^\infty(\Omega)$ norm.}
    \label{Fig:H14}
\end{figure}

We apply standard finite differences with $n=200$ grid points to discretize the problem in each spatial direction.
The reference solution is computed by the ERKC-C method with Gauss collocation points ($s=3$) using constant step size $h=2^{-14}$. The temporal orders of convergence of ERKC-I and ERKC-C methods in the $L^\infty(\Omega)$ norm at final time $T=3$ are presented in Figure~\ref{Fig:H14}. The observed orders evidently are in line with our theoretical analysis.

\begin{example}\label{Exa:mERKC-I}
In the last example, we test the ERKC-I method and its modified version on the semilinear parabolic problem subject to periodic boundary conditions 
\begin{equation}\label{Eqn:mERKC}
\frac{\partial u}{\partial t} (t,x) -\frac{\partial^2 u}{\partial x^2}(t,x) = \frac{1}{1+u(t,x)^2}+\frac{2000}{1+u(t^2-1,x)^2}+\Phi(t,x).
\end{equation}
We consider this problem for $x\in[0,1]$ and $t\in [0,1.4]$. The source function $\Phi$ is determined by the exact solution of the problem $u(t,x)= \Psi(t) \sin (2\pi x)$, where $\Psi(t)$ is given by
$$
\Psi(t) = 
\left\{\begin{aligned}
& \mathrm{e}^{-t}, && t\in[-1,0], \\
& 1+t\mathrm{e}^{2t}, && t\in(0,1], \\
& (1+\mathrm{e}^{2})+3\mathrm{e}^{2}(t-1)+(t-1)^2\mathrm{e}^{3t}, && t \in(1,1.4].
\end{aligned}\right.
$$  
\end{example}

\begin{figure}[H]
    \centering
    \includegraphics[width=0.49\textwidth]{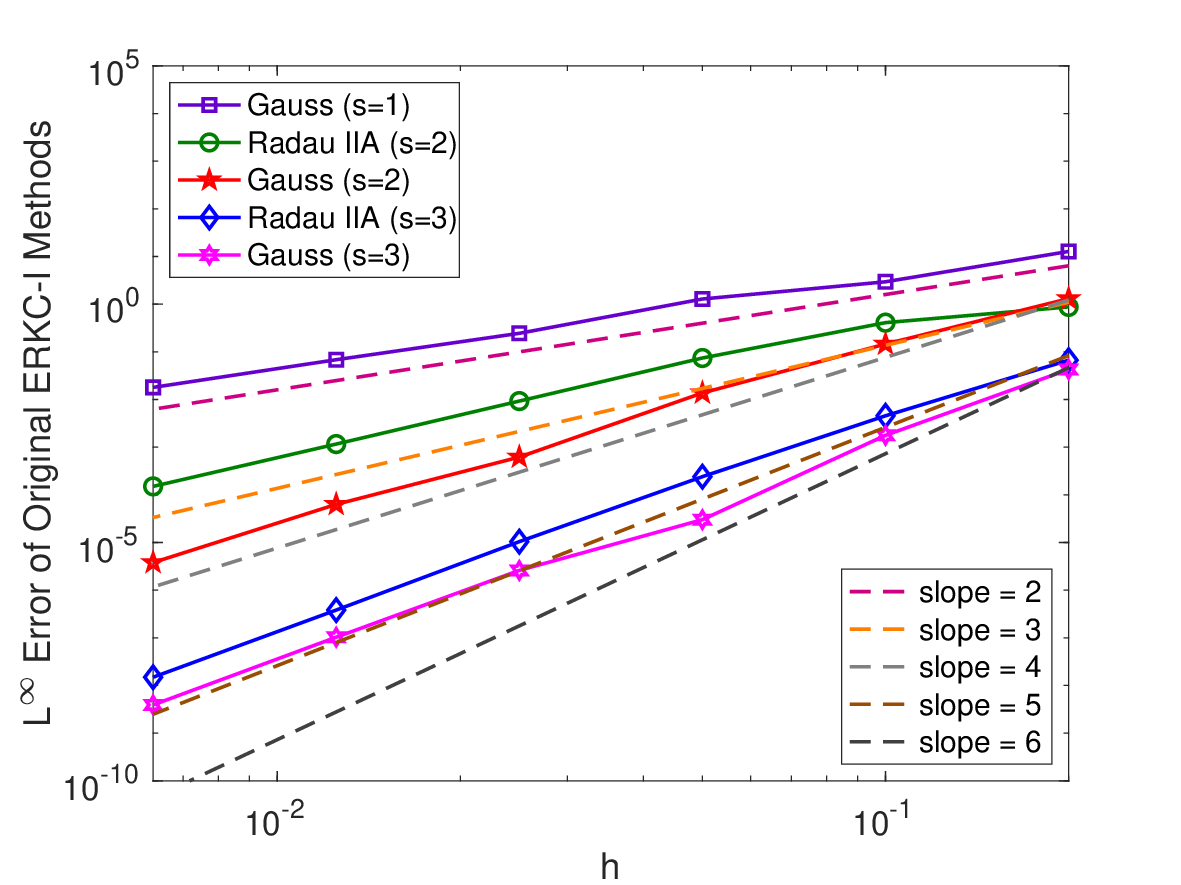}
    \includegraphics[width=0.49\textwidth]{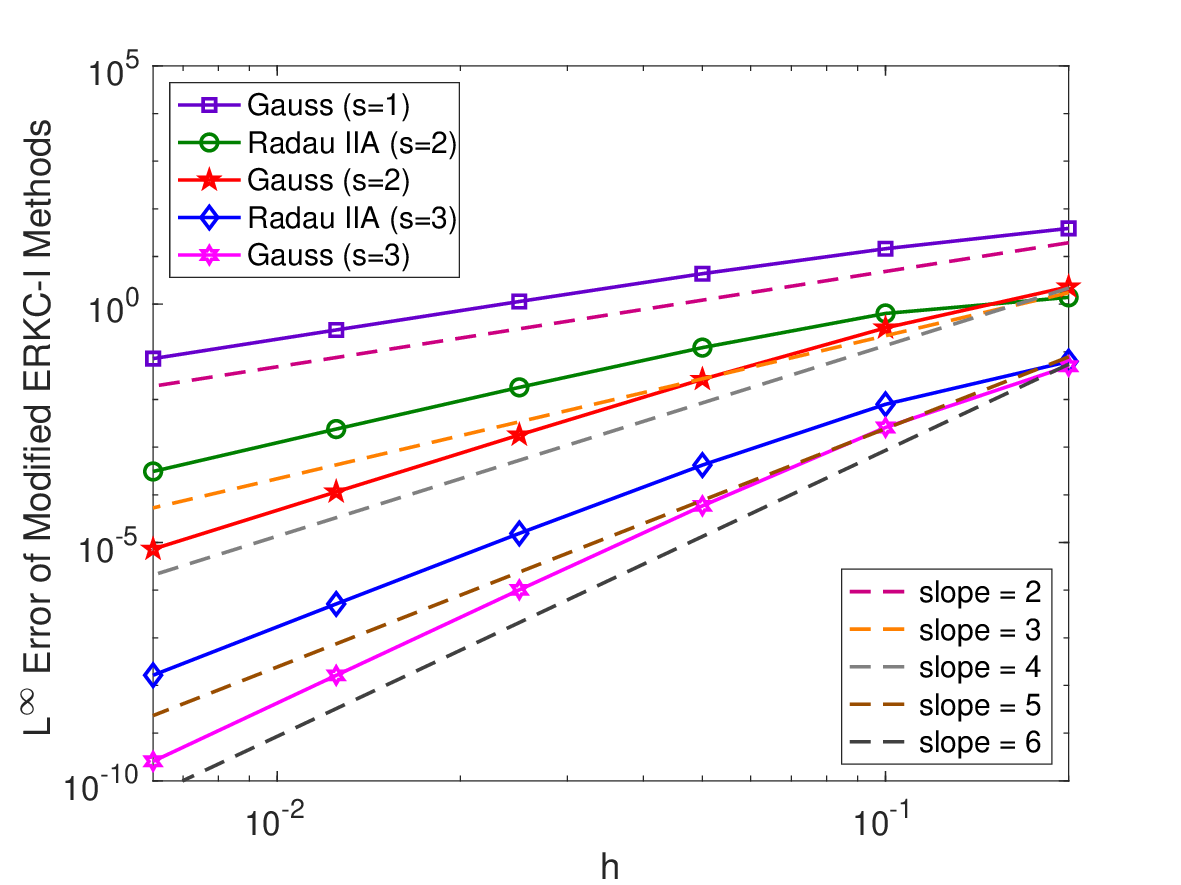}
    \caption{The convergence rates of ERKC-I methods (in the left panel) and modified ERKC-I methods (in the right panel) for \eqref{Eqn:mERKC}. The errors are measured at $T=1.4$ in the $L^\infty(\Omega)$ norm.}
    \label{Fig:nd}
\end{figure}

We apply a standard pseudospectral method with $n=200$ grid points to discretize problem \eqref{Eqn:mERKC} in space. The convergence rates of the ERKC-I methods and their modified versions at $T = 1.4$ in the $L^\infty(\Omega)$ norm, respectively, are presented in Figure \ref{Fig:nd}. It is observed that the modified ERKC-I methods achieve full order accuracy, while the original ERKC-I methods fail when $s=3$.  However, ERKC-I methods generally exhibit superior accuracy compared to their modified counterpart, except for $s = 3$ with Gauss collocation points.

\section{Conclusion}
In this paper, we extended the ERKC methods to semilinear parabolic problems with time-dependent delay. Two classes of ERKC methods, ERKC-I and ERKC-C methods, were proposed. 
Considering the typically nonsmooth nature of solutions to delay differential equations, we employed constrained meshes to ensure that potential discontinuity points were captured.
We proved that the continuous numerical solutions converge to the true solutions with order $s$ at least. Order $s+1$ can be achieved provided that the underlying quadrature rule is of order $s+1$. 
A modified ERKC-I method was constructed, which can further achieve $s+1+\beta$, provided that the underlying quadrature rule is of order $s+2$ and some slightly stronger assumptions regarding spatial regularity hold. 
A potential future direction could be the extension of ERKC method to semilinear parabolic problems with vanishing delay.

\section*{Acknowlegments}
Qiumei Huang is supported by the National Natural Science Foundation of China (No.~12371385). 
Gangfan Zhong is supported by the China Scholarship Council (CSC) joint Ph.D. student
scholarship (Grant 202406540082).

\end{document}